\documentclass[final, 3p, times]{elsarticle}

\usepackage{graphicx, hyperref, xspace}
\usepackage{amsmath, amsfonts, amssymb, xcolor, comment, bm}

\newcommand{\kiops}{\texttt{KIOPS}\xspace}

\newcommand{\epic}{\texttt{EPIC}\xspace}

\usepackage{tikz}
\usetikzlibrary{calc,shapes}

\newcommand{\tikzmark}[1]{\tikz[overlay,remember picture] \node (#1) {};}

\usepackage[caption=false]{subfig}
\hypersetup{colorlinks,
            linkcolor = {red!50!black},
            citecolor = {red!10!black},
            urlcolor  = {green!10!black}}

\newcommand*{\colorboxed}{}
\def\colorboxed#1#{%
  \colorboxedAux{#1}%
}
\newcommand*{\colorboxedAux}[3]{%
  \begingroup
    \colorlet{cb@saved}{.}%
    \color#1{#2}%
    \boxed{%
      \color{cb@saved}%
      #3%
    }%
  \endgroup
}

\journal{Computer Physics Communications}

\begin{document}

\begin{frontmatter}



\title{A comparison of Leja-- and Krylov--based iterative schemes \\ for Exponential Integrators}


\author[inst1]{Pranab J. Deka} \corref{lod1}
\ead{pranab.deka@uibk.ac.at}
\cortext[lod1]{Corresponding author}

\author[inst2]{Mayya Tokman}
\author[inst1]{Lukas Einkemmer}

\affiliation[inst1]{organization = {Department of Mathematics, University of Innsbruck},
            city = {Innsbruck},
            postcode = {A-6020}, 
            country = {Austria}}
            
\affiliation[inst2]{organization = {School of Natural Sciences, University of California},
            city = {Merced},
            postcode = {CA 95343}, 
            country = {USA}}

\begin{abstract}
Krylov-based algorithms have long been preferred to compute the matrix exponential and exponential-like functions appearing in exponential integrators. Of late, direct polynomial interpolation of the action of these exponential-like functions have been shown to be competitive with the Krylov methods. We analyse the performance of the state-of-the-art Krylov algorithm, \kiops, and the method of polynomial interpolation at Leja points for a number of exponential integrators for various test problems and with varying amounts of stiffness. Additionally, we investigate the performance of an iterative scheme that combines both the \kiops and Leja approach, named LeKry, that shows substantial improvements over both the Leja- and Krylov-based methods for certain exponential integrators. Whilst we do manage to single out a favoured iterative scheme for each of the exponential integrators that we consider in this study, we do not find any conclusive evidence for preferring either \kiops or Leja for different classes of exponential integrators. We are unable to identify a superior exponential integrator, one that performs better than all others, for most, if not all of the problems under consideration. We, however, do find that the performance significantly depends on the interplay between the iterative scheme and the specific exponential integrator under consideration.
\end{abstract}




\end{frontmatter}



\section{Introduction}

Exponential integrators have shown remarkable promise and have been receiving increased attention from the scientific community for solving time-dependent partial differential equations (PDEs) of the form  
\begin{equation}
    \frac{\partial u}{\partial t} = f(u(t)), \qquad u(t_0) = u_0.
    \label{eq:nl_pde}
\end{equation}
In addition to possessing excellent stability conditions, they tend to be highly accurate whilst being able to take substantially larger step sizes than implicit integrators. They are of particular importance in systems of PDE(s) that involve different characteristic time scales and the fastest time scale of the system is much shorter than the slowest one, in essence, stiff systems. An excellent review on exponential integrators have been presented by \citet{Ostermann10}.

Of particular interest are the Exponential Propagation Iterative Runge--Kutta (EPIRK \cite{Tokman06, Tokman12}) and Exponential Rosenbrock (EXPRB \cite{Caliari09, Ostermann10}) methods that linearise the underlying PDE(s) at every time step: the linear term is solved exactly (in time) and the nonlinear term is approximated with some explicit methods. As an example, Eq. \eqref{eq:nl_pde} can be linearised as
\begin{equation}
    \frac{\partial u}{\partial t} = \mathcal{J}(u^n) \, u^n + \mathcal{F}(u^n), \nonumber
\end{equation}
where $\mathcal{J}(u^n)$ is the Jacobian of the nonlinear function $f(u^n)$ at the $n^\mathrm{th}$ time step and $\mathcal{F}(u^n) = f(u^n) - \mathcal{J}(u^n) \, u^n$ is the nonlinear remainder. Worth noting is that a multi-stage exponential integrator can have several internal stages. For some internal stage, say $k$, we define the difference of these nonlinear remainders as $\mathcal{R}(k) = \mathcal{F}(k) - \mathcal{F}(u)$.

Exponential integrators, as the name suggests, involve computing the exponential of a matrix or exponential-like functions of the matrices related to the system under consideration. These exponential-like functions, or the $\varphi$ functions are given by $\varphi_{l + 1}(z) = \frac{1}{z} \left(\varphi_l(z) - \frac{1}{l!} \right), l \geq 1$ and $\varphi_0(z) = e^z$. Several methods for computing the matrix exponential have been reviewed by Moler and Loan \cite{Moler78, Moler03}. Although Pad\'e approximation, scaling and squaring, or diagonalising the matrix are excellent choices to compute these exponential-like functions for small matrices, they become prohibitively expensive for large-scale stiff systems. In case of large systems of matrices, Krylov methods have been proposed \cite{Vorst87} and are commonly used to compute the exponential-like functions that appear in exponential integrators. In addition to the remarkable ability of the Krylov-based methods to effectively handle systems of large matrices, they require no prior information on the spectrum of the matrix. The general idea of Krylov-based algorithms is to project the exponential of a large matrix onto a small Krylov subspace and compute the exponential of the smaller matrix using computationally cheap methods. Another class of methods that have been shown to be competitive with the Krylov methods is the method of polynomial interpolation at Leja points. Whilst this method does need some estimate of the spectrum of the underlying matrix, they are relatively simpler and easier to implement than the Krylov-based ones. The key advantage of the Leja method is that one needs to store only one input and one output vector, in essence,  their memory requirements are similar to that of explicit schemes (as opposed to Krylov methods where one needs to potentially store a relatively large number of Krylov vectors). This makes Leja methods very well suited for graphic processing units (GPUs). This is particularly important as most of the fastest supercomputers are now based on this technology. Thus, to obtain good performance on such systems a method (such as Leja interpolation) that requires as little communication as possible has to be used. The need for the evaluation of inner products in Krylov-based methods can result in performance loss when scaling such methods to massively parallel GPU systems or large supercomputers.

In \citet{Deka22a}, we studied the performance of the Leja and Krylov-based methods for the set of equations of magnetohydrodynamics (MHD). We found that Leja had better performance than the Krylov method for the EXPRB43 \cite{Caliari09, Hochbruck09, Ostermann10} integrator whereas Krylov had surpassed Leja for the ``Krylov-friendly'' EPIRK5P1 \cite{Tokman12} integrator. Based on these results, we hypothesised that since the EPIRK class of integrators are designed to be ``Krylov-friendly'', Krylov-based methods are expected to be favourable for this class. In this manuscript, we aim to test this hypothesis for a wide range of test problems for various EPIRK and EXPRB integrators using the Leja and Krylov-based methods. In Sec. \ref{sec:iterative_schemes}, we provide a gist of the Leja and Krylov methods as well as our variable step size implementation for improved computational performance, which is crucial for the Leja-based methods. Let us clearly state that we use the state-of the-art \kiops algorithm \cite{Gaudreault18} for this work, which is different from the elementary Krylov-based method that we used in the previous study \cite{Deka22a}. We describe the test problems in Sec. \ref{sec:problems} and report further developments in the Leja algorithm, in Sec. \ref{sec:vertical}, that result in considerable speedups without any loss in accuracy. A comparison of the performance of Leja, \kiops, and an iterative scheme based on the combination of the two for the aforementioned range of problems is presented Sec. \ref{sec:lekry}. We compare the performance of the different exponential integrators in Sec. \ref{sec:compare_expint} and finally conclude in Sec. \ref{sec:conclude}.


\section{Iterative schemes and step size control}
\label{sec:iterative_schemes}

In this study, we use Exponential Propagation Integrators Collection (\epic), a publicly available package (\url{https://faculty.ucmerced.edu/mtokman/#software}) to compare the computational performances of the Leja method with the widely used Krylov-subspace algorithm. We provide a gist of these algorithms in the following subsection and refer to the respective original papers, cited accordingly, for a detailed description. All simulations in this work have been conducted in \texttt{MATLAB} and proper care has been taken to avoid superfluous computational expenses for both the Leja- and Krylov-based algorithms. The major computational cost in either of the  algorithms is dictated by cost of the matrix-vector products, or in the case of matrix-free implementation, as in this work, the computation of the right-hand-side (RHS) functions.


\subsection{Krylov subspace algorithm}

\citet{Niesen09} proposed the \texttt{phipm} algorithm for computing the exponential of a matrix and the $\varphi$ functions appearing in exponential integrators. \texttt{phipm} reduces the computation of a large matrix, usually arising from the discretisation of a PDE, by projecting the matrix onto a small Krylov subspace. This is usually done using Arnoldi iterations \citep{Arnoldi1951}, thereby forming the  Hessenberg matrix. The task of computing the exponential-like functions of a large matrix is reduced to computing the exponential-like functions of the smaller Hessenberg matrix. Additionally, \texttt{phipm} computes the exponential of a larger augmented matrix over computing multiple small $\varphi_l(z)$ functions \cite{Higham10}, as demonstrated in the following theorem: \\
\textbf{Theorem}: Given $A \in \mathbb{R}^{N \times N}$, $B = [b_p, \hdots, b_2, b_1] \in \mathbb{R}^{N \times p}$, $K = \begin{bmatrix} 0 & I_{p-1} \\ 0 & 0 \end{bmatrix} \in \mathbb{R}^{p \times p}$, $v = [b_0, e_p] \in \mathbb{R}^{N+p}$, we define an augmented matrix, $\Tilde{A} = \begin{bmatrix} A & B \\ 0 & K \end{bmatrix} \in \mathbb{R}^{(N + p) \times (N + p)}$ and $w = \exp(\Tilde{A} \Delta t) v$. The first $N$ elements of $w$ are given by $w(1:N) = \sum_{j = 0}^{p} \varphi_j(A \Delta t) b_j \Delta t^j$. \\
The proof of this theorem can be found in the aforementioned papers. Writing the linear combination of $\varphi_l(z)$ functions as an augmented matrix exponential results in a few advantages: First, the number of exponential-like functions, needed to be evaluated, are reduced. Second, one can divide the step size into multiple small substeps that is expected to result in overall computational savings. Third, since any stage of an exponential integrator can be written as a linear combination of $\varphi_l$ functions, one can easily generalise the algorithm to compute any linear combination of `n' $\varphi_l$ functions. An error estimate of the Krylov subspace algorithm, substepping of the time step, and adaptive modification of the dimension of the Krylov subspace are used to reduce the error incurred and to prevent the Krylov basis from growing too large as this can severely hamper the computational performance.

\citet{Gaudreault18} proposed the \kiops algorithm that build on the ideas of \texttt{phipm}. Similar to \texttt{phipm}, \kiops computes an approximation of the exponential on the Krylov subspace and considers substepping to compute the exponential of the augmented matrix. In addition, they improved on the adaptivity procedure of selecting the step sizes for substepping and the size of the Krylov subspace that results in an overall improvement in the efficiency of the algorithm. Incomplete orthogonalisation was proposed to combat the computational expense of the Arnoldi procedure for the Krylov projections. This reduces the number of inner products needed per Krylov iteration to as few as two. The use of incomplete orthogonalisation significantly alleviates the challenges associated with parallelisation of the computation of inner products. It is to be noted that Krylov-based methods do not require any prior estimates of the spectrum of the underlying matrix.

Whilst the main computational cost in the \kiops algorithm arises from the need to compute the matrix-vector products or the RHS function, computing the exponential of the Hessenberg matrix and step size rejections do play a moderate role in the overall expenses of the algorithm.


\subsection{Leja polynomial interpolation}

The Leja polynomial interpolation method for exponential integrators \cite{Caliari04, Bergamaschi06, Caliari07b, Caliari09} involves interpolating the $\varphi_l(z)b_k$ on a set of predetermined Leja points \cite{Leja1957, Baglama98}. One requires a rough estimate of the spectrum of the underlying Jacobian matrix which can be obtained using power iterations. It is to be noted that if the solution is evolved in time by relatively small step sizes, the solution, and consequently the spectrum is not expected to changes drastically. Taking this into account, one can compute an approximation to the spectrum of the Jacobian every `n' time steps, where the spectrum is subject to a safety factor to ensure an overestimation of the largest eigenvalue. Say, the largest and the smallest real eigenvalue, in magnitude, are $\alpha$ and $\beta$ respectively. Mathematically, this can be written as $\alpha \leq \text{Re} \, \sigma(z) \leq \beta \leq 0$. Now, we scale and shift the spectrum of the matrix on the set of pre-computed Leja points ($\xi$) in the arbitrary domain [-2, 2]. Then, we compute the coefficients of the function to be interpolated, i.e. $\varphi(\Delta t (c + \gamma \xi))$, where $c$ is the midpoint of the spectrum and $\gamma$ is one-forth the distance of the distance between $\alpha$ and $\beta$. The $m^\mathrm{th}$ term of the polynomial, $p_m(z)$, is given by
\begin{align*}
    p_{m+1}(z) & = p_m(z) + d_{m+1} \, y_{m+1}(z), \\
    y_{m+1}(z) & = y_m(z) \times \left(\frac{z - c}{\gamma} - \xi_{m} \right),
\end{align*}
This iterative procedure, i.e. computation of the RHS function ($z$), constitutes the most expensive part of the Leja interpolation method. Each of the stages of an exponential integrator is interpolated as a polynomial. We use a termination criterion that stops the interpolation procedure if the error incurred ($|d_m| \, \|y_m\|$) is less than a certain user-defined tolerance. The convergence of the augmented matrix for Leja interpolation is remarkably slow. Furthermore, considering an augmented matrix (and substepping) for the Leja interpolation case, would result in having to compute and store additional coefficients for the polynomials (divided differences) and multiple intermediate vectors. This leads to a severe deterioration in the computational performance. We show this in \ref{app:linear_Leja_phi}. This is why we compute the interpolation of the individual $\varphi_l$ functions for this iterative scheme.

The number of Leja points needed for convergence depends on, among other parameters, the step size. Larger step sizes require more Leja points for the polynomial to converge. This is impractical as the coefficients, determined using the divided differences algorithm, quickly goes down to machine precision. The algorithm, thus, becomes prone to round-off errors. The problem becomes even more severe when the norm of the function to be interpolated is large, i.e., $\varphi_1(\mathcal{J}(u)) f(u) \Delta t$. The convergence of this function constitutes the bottleneck for convergence of an iterative method at any given step. Note that substepping can only be used for the computation of the matrix exponential and not for $\varphi_l$ functions. Thus, it becomes necessary to choose a relatively small step size, wherever possible, to avoid step size rejections. 


\subsection{Automatic step size control}

One of the crucial aspects of automatic time step control is the selection of an appropriate step size controller. Traditional step size controllers, such as the PI controller, consider the largest possible step size subject to the prescribed accuracy constraints. The step size at the $(n+1)^\mathrm{th}$ time step ($\Delta t^{n+1}$) is given by
\[
    \Delta t^{n+1} = \Delta t^n \times \left(\frac{\mathrm{tol}}{e^n} \right)^{1/(p + 1)},
\]
where $\Delta t^n$ and $e^n$ are the step size and an estimate of the error incurred at the $n^\mathrm{th}$ time step, respectively, tol is the user-defined tolerance, and $p$ is the order of convergence of the integrator. The major drawback of such a controller is the implicit assumption of the independence of the computational cost on the step size. This is true only for explicit integrators and direct solvers. In the case of iterative solvers (for implicit and exponential integrators), the computational cost incurred is a function of the step size: large step sizes need more iterations to converge to the prescribed tolerance as compared to smaller step sizes. The dependence of the cost incurred on the step size is (highly) nonlinear. \citet{Einkemmer18} developed a step size controller that minimises the computational cost by choosing a smaller step size whenever possible. This step size controller, or the cost controller from here on, has been optimised to incur minimum possible computational cost by choosing a smaller step size, than the one chosen by the traditional controller. The cost controller was first proposed for implicit integrators \cite{Einkemmer18} and was later shown to be even more beneficial for exponential integrators \cite{Deka22a, Deka22b}. This is because exponential integrators, in principle, are able to take much larger step sizes than implicit ones. Consequently, opting for smaller step sizes than the maximum determined from the error estimate can result in drastic computational savings. Additionally, choosing smaller step sizes is likely to result in more accurate solutions. Therefore, one can expect better results in terms of accuracy, for this controller. For details on the optimisation procedure and the derivation of the equations, we refer the reader to the original work \cite{Einkemmer18}. The equation for the cost controller reads
\[
    \Delta t^{n+1} = \Delta t^n \times
    \begin{cases}
    
    \lambda                                                 & \text{if $1 \leq s < \lambda$}, \\
    \delta                                                  & \text{if $ \delta \leq s < 1$}, \\
    s:=\exp({-\alpha \; \mathrm{tanh}(\beta \Delta)})       & \text{otherwise},
    \end{cases}
\]
where $\Delta = (\mathrm{ln}\, c^n - \mathrm{ln}\, c^{n - 1})/(\mathrm{ln}\, \Delta t^n - \mathrm{ln}\, \Delta t^{n - 1})$, $c^n$ is the computational cost incurred at the $n^\mathrm{th}$ time step and $\alpha = 0.65241444$, $\beta= 0.26862269$, $\lambda = 1.37412002$, and $\delta = 0.64446017$. The cost controller only reduces the computational cost and does not take into account the accuracy constraints. As such, one needs to choose the minimum of the step sizes yielded by the cost and the traditional controller, i.e. $\Delta t = \text{min}(\Delta t_\mathrm{traditional}, \Delta t_\mathrm{cost})$. The cost controller has been used, to reduce the computational runtimes, for all simulations conducted in this work. 


\section{Test Problems}
\label{sec:problems}

We used the following problems to demonstrate the enhanced `Vertical' implementation in Leja (Sec. \ref{sec:vertical}), compare the efficacy of the \texttt{KIOPS} algorithm with that of the Leja method (Sec. \ref{sec:lekry}), and compare the relevant exponential integrators (Sec. \ref{sec:compare_expint}). The test examples have been drawn from Loffeld \& Tokman \cite{Tokman13}, Rainwater \& Tokman \cite{Tokman16}, and Gaudreault et al. \cite{Gaudreault18}. In all of the following problems, $\nabla$ and $\nabla^2$ are discretized using the second-order centered finite difference scheme.


\paragraph{Advection-–diffusion–-reaction (2D)}

The advection--diffusion--reaction (ADR) equation reads \citep{Caliari09}
\begin{equation*}
    \frac{\partial u}{\partial t} = \alpha \nabla^2 u + \beta \nabla u + \gamma \, u \, (1 - u) \, \left(u - \frac{1}{2} \right), \quad x, y \in [0, 1], \quad t \in [0, 0.01] ,
\end{equation*}
with homogeneous Neumann boundary conditions. The initial condition is chosen to be
\begin{equation*}
	u(x, y, t = 0) =  256(xy \, (1 - x)(1 - y ))^2 + 0.3.
\end{equation*}
We consider two different values of the diffusion coefficient: (a) $\alpha = 0.1$, $\beta = -10$, and $\gamma = 100$, and (b) $\alpha = 0.01$, $\beta = -10$, and $\gamma = 100$.


\paragraph{Allen--Cahn (2D)}

The Allen--Cahn equation can be written as \cite{Bates09}
\begin{equation*}
    \frac{\partial u}{\partial t} = \alpha \nabla^2 u + u - u^3, \quad x, y \in [-1, 1], \quad t \in [0, 1.0],
\end{equation*}
with no-ﬂow boundary conditions. We choose three different values of $\alpha = 10^{-1}, 10^{-2}, 10^{-3}$ to investigate different stiffness regimes. The initial condition is given by
\begin{equation*}
	u(x, y, t = 0) =  0.1 + 0.1 \cos(2 \pi x) \cos(2 \pi y).
\end{equation*}


\paragraph{Brusselator (2D)}

The two-dimensional Brusselator equation is \cite{Lefever71, HairerII}
\begin{align*}
    \frac{\partial u}{\partial t} & = \alpha \nabla^2 u + uv^2 - 4u + 1, \quad x, y \in [0, 1], \quad t \in [0, 1.0],  \\
    \frac{\partial v}{\partial t} & = \alpha \nabla^2 u - u^2v + 3u, 
\end{align*}
with homogeneous Neumann boundary conditions. We consider three different values of diffusion coefficients: $\alpha = 10^{-1}, 10^{-2}, 10^{-3}$. Initial conditions are chosen to be
\begin{align*}
    u(x, y, t = 0) & =  2 + 0.25 y, \\
    v(x, y, t = 0) & =  1 + 0.8 x.
\end{align*}


\paragraph{Gray--Scott (2D)}

The Gray--Scott equations read \cite{Gray84}
\begin{align*}
    \frac{\partial u}{\partial t} & = \alpha_1 \nabla^2 u - uv^2 + a(1 - u), \quad x, y \in [0, 1], \quad t \in [0, 0.1], \\
    \frac{\partial v}{\partial t} & = \alpha_2 \nabla^2 v + uv^2 - (a + b)v
\end{align*}
with periodic boundary conditions, $a = 0.04$, $b = 0.06$, and $\alpha_1, \alpha_2 = 10^{-1}, 10^{-2}, 10^{-3}$. The initial conditions are
\begin{align*}
    u(x, y, t = 0) & =  1 - \exp{\left[-150 \left(\left(x - \frac{1}{2}\right)^2 + \left(y - \frac{1}{2}\right)^2\right)\right]}\\
    v(x, y, t = 0) & =  \exp{\left[-150 \left(\left(x - \frac{1}{2}\right)^2 + 2\left(y - \frac{1}{2}\right)^2\right)\right]}.
\end{align*}


\paragraph{Semilinear (1D)}

The one-dimensional semilinear parabolic equation \cite{Hochbruck05} is
\begin{equation*}
    \frac{\partial u}{\partial t} = \frac{\partial^2 u}{\partial x^2} + \int_0^1 u \, dx + \Phi (x, t),
\end{equation*}
where we consider homogeneous Dirichlet boundary conditions for $x \in [0, 1]$ and the simulation time, $t_f = 1$. $\Phi(x, t)$ is chosen such that $u(x, t) = x(1 - x) e^t$ is the exact solution. 


\section{Vertical and Mixed Implementation in Leja and Krylov}
\label{sec:vertical}

Rainwater \& Tokman \cite{Tokman16} introduced the ``vertical exponential Krylov'', ``horizontal exponential adaptive Krylov'', and ``mixed exponential adaptive Krylov'' methods as means of further reducing the computational cost whilst implementing a Krylov-based exponential integrator. The vertical implementation is based on the idea of computing the columns of each stage that share the same vector $b_k$, simultaneosuly. For instance, in EPIRK4s3 (Eq. \eqref{eq:epirk4s3}), the functions $\varphi_1\left(\frac{1}{8} \mathcal{J}(u^n)  \Delta t \right)$, $\varphi_1\left(\frac{1}{9} \mathcal{J}(u^n)  \Delta t \right)$, and $\varphi_1\left(\mathcal{J}(u^n)  \Delta t \right)$ are applied to the same vector $f(u^n)\, \Delta t$, and as such these three projections can be computed at once. The computational cost of this is equivalent to the cost of computing only $\varphi_1\left(\mathcal{J}(u^n) \Delta t \right)\,f(u^n)\, \Delta t$; $\varphi_1\left(\frac{1}{8} \mathcal{J}(u^n)  \Delta t \right)\,f(u^n)\, \Delta t$ and $\varphi_1\left(\frac{1}{9} \mathcal{J}(u^n)  \Delta t \right)\,f(u^n)\, \Delta t$ can be obtained with a negligible overhead.

\begin{figure}[ht]
    \centering
	\includegraphics[width = \textwidth]{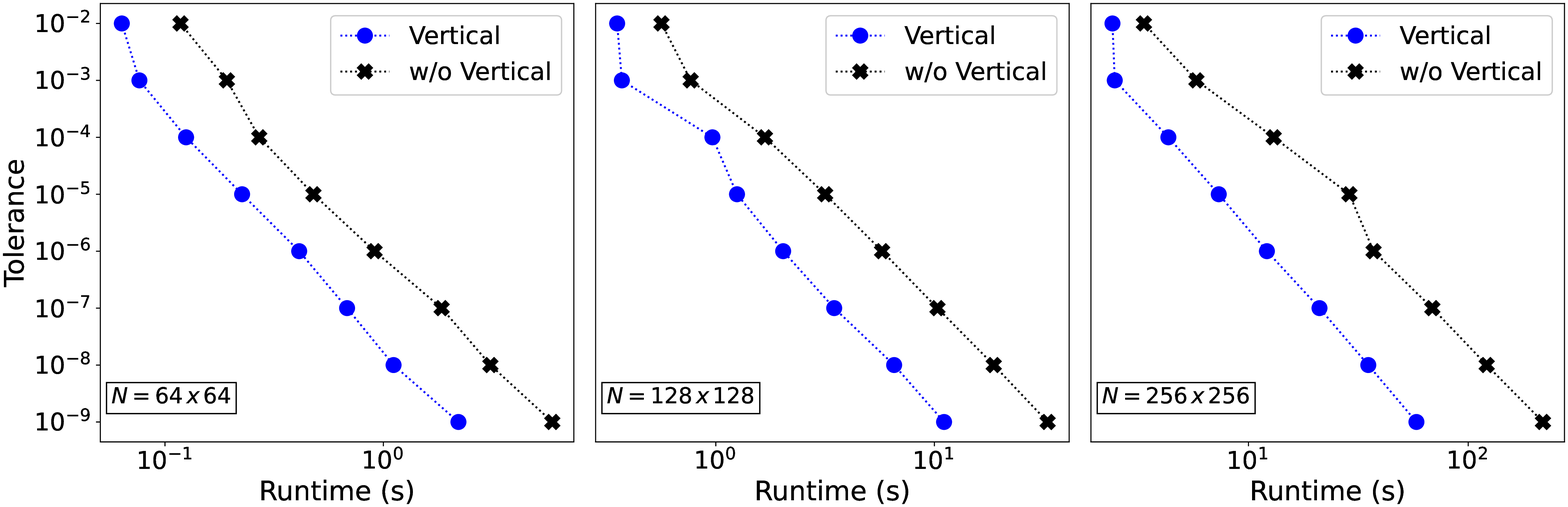} \\
	\includegraphics[width = \textwidth]{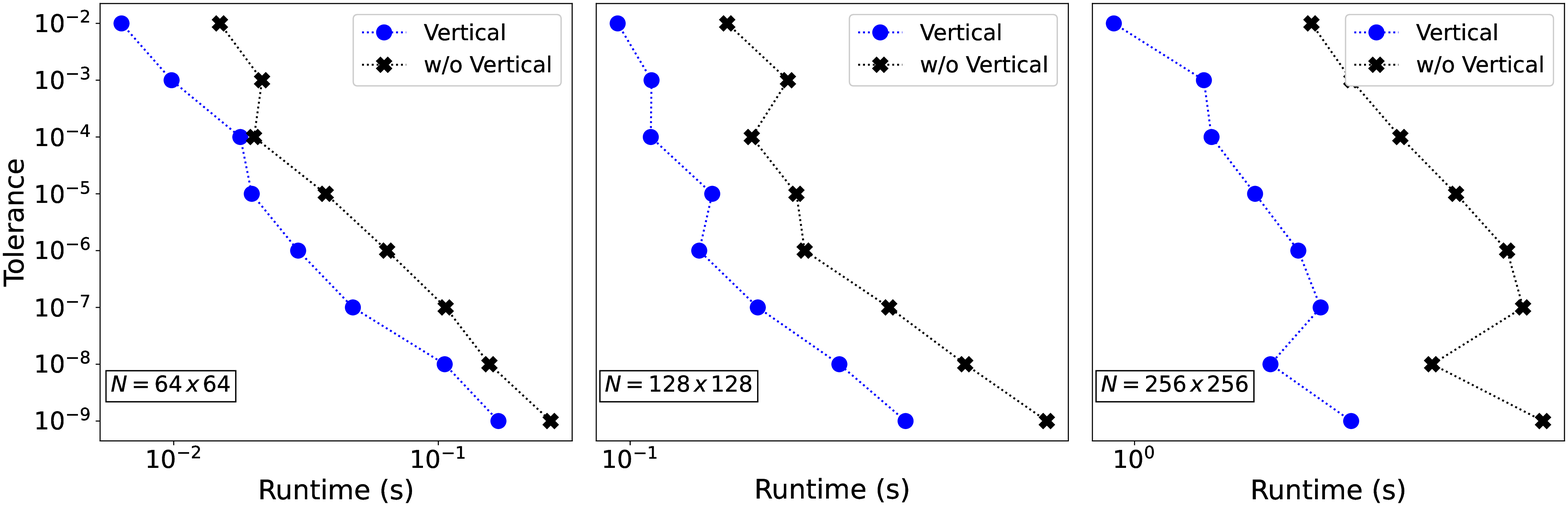}
    \caption{Comparison of the runtimes with (blue circles) and without (black crosses) the vertical implementation in the Leja method with EPIRK4s3 for the Brusselator ($\alpha = 10^{-3}$, top panel) and the Allen--Cahn ($\alpha = 10^{-2}$, bottom panel) equations. Clearly, the vertical implementation results in a significant amount of computational savings.}
    \label{fig:compare_Leja_nonvertical}
\end{figure}

For the Leja method, we propose a similar vertical approach. However, there is one major difference: the integrator coefficients, here 1/8, 1/9, and 1, determine the coefficients of the polynomial being interpolated. Therefore, one needs to construct three separate vectors for these three Leja interpolations. Nevertheless, one can still obtain a substantial amount of computational savings by going through the iteration procedure, i.e. computing the numerical Jacobian and the re-scaling and re-shifting only once (instead of thrice, as would be the case if these three interpolations were done separately). This is a novel approach for the Leja interpolation scheme and we highlight in Fig. \ref{fig:compare_Leja_nonvertical} that this leads to significant computational savings.

\begin{figure}[ht]
    \centering
	\includegraphics[width = \textwidth]{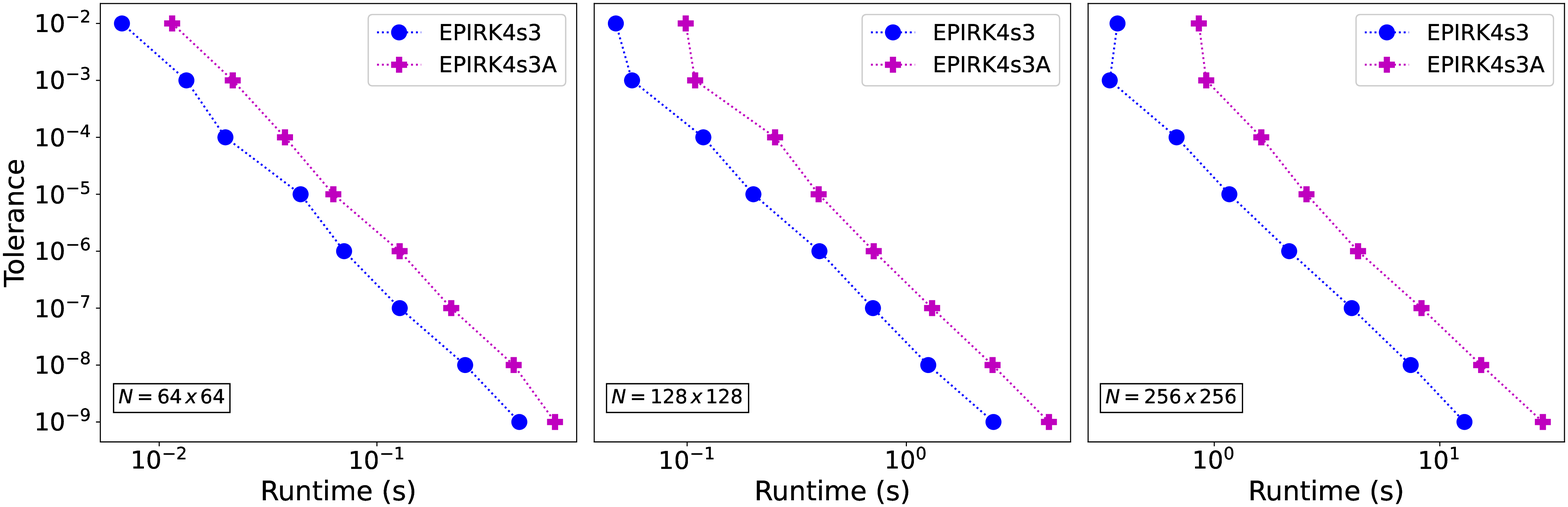} \\
	\includegraphics[width = \textwidth]{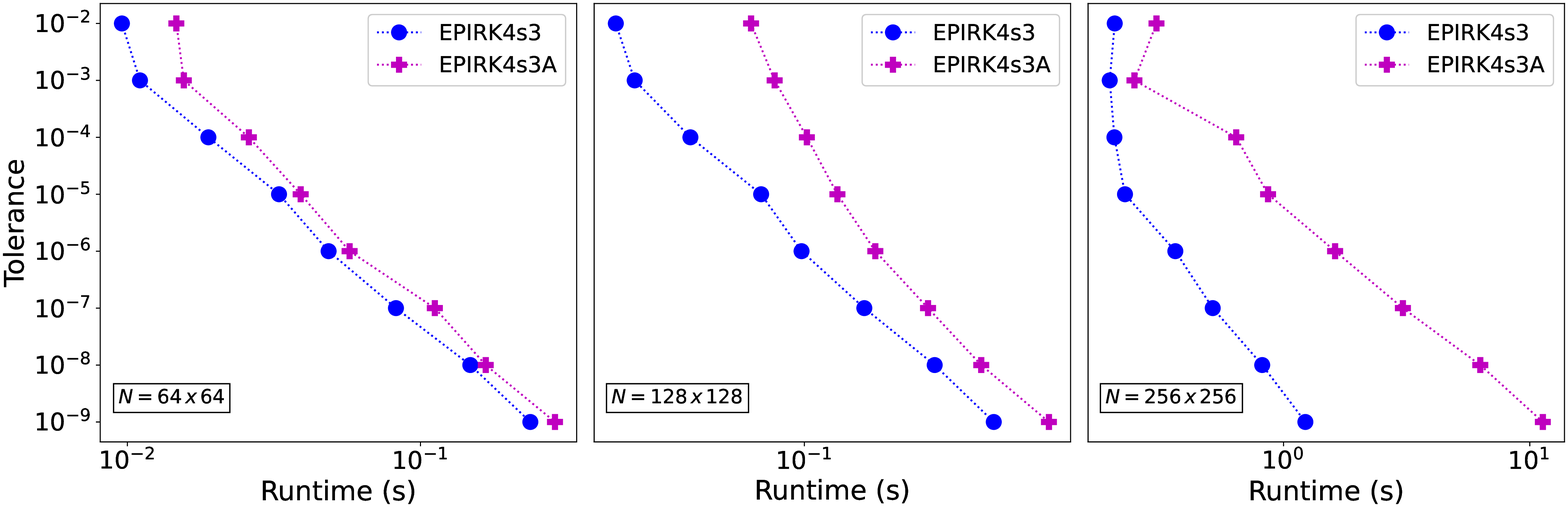}
    \caption{Comparison of the runtimes of the internal stages (in vertical) of EPIRK4s3 and EPIRK4s3A, with the Leja scheme, for the Brusselator ($\alpha = 10^{-3}$, top panel) and the Allen--Cahn ($\alpha = 10^{-2}$, bottom panel) equations. It can be seen that smaller integrator coefficients (in EPIRK4s3) result in faster convergence.}
    \label{fig:compare_LeKry_vertical}
\end{figure}

The horizontal implementation is intended to compute all terms, at any given stage, with one Krylov evaluation. The mixed implementation is, as the name suggests, a mixture of vertical and horizontal implementations. Typically, the internal stages are computed vertically whereas the final stage is computed horizontally. The vertical implementation of the internal stages, in a mixed implementation, can reduce the cost of computing the internal stages by a substantial amount. This is due to the fact an integrator coefficient, say $l$ ($< 1$), effectively scales down the Jacobian to a factor of $l$. Smaller values of $l$ insinuate faster convergence and reduced cost. The reduction in the cost scales nonlinearly with the integrator coefficient. The largest integrator coefficients of the internal stages of EPIRK4s3 and EPIRK4s3A (Eq. \eqref{eq:epirk4s3a}) are 1/8 and 2/3, respectively, and as such the internal stages, in the vertical implementation, of EPIRK4s3 are expected to converge faster than EPIRK4s3A (Fig. \ref{fig:compare_LeKry_vertical}). These results are in agreement with \citet{Tokman17b} where they compared the performance of these two integrators for constant step sizes. In this comparison, we are assuming that since both EPIRK4s3 and EPIRK4s3A have an identical structure of equations and they are of the same order, the effects of variable step sizes is negligible. One has to take care whilst comparing the runtimes of the internal stages of different integrators in vertical - if two integrators have different different orders of convergence (owing to ability of the higher-order integrator to take larger step sizes) or different number of internal stages, this comparison would be rendered invalid. Rainwater \& Tokman \cite{Tokman16} developed several integrators tailored to these vertical, horizontal, and mixed implementations. Let us emphasise that any integrator could be implemented in one of these three ways, however, the efficiency of these are likely to differ based on the integrator coefficients at the various stages.


\section{Comparing and combining \kiops and Leja for Exponential Integrators}
\label{sec:lekry}

In this section, we compare the performance of the \kiops algorithm with that of Leja, and propose an efficient way of combining these two algorithms for a given integrator. It is to be noted that the most efficient implementation for an integrator may differ for constant and variable step size implementations. Since we are primarily interested in variable step sizes, we focus on embedded exponential integrators, in essence, exponential integrators with an error estimator, that shares the same internal stages, embedded in them. We focus on the EPIRK4s3, EPIRK4s3A, EPIRK5P1, EXPRB43, and EXPRB53s3 integrators, the details of which can be found in the following subsections.


\newcommand{\Drawa}[3]{%
  \begin{tikzpicture}[overlay,remember picture]
    \draw[#3] ([yshift = 15pt, xshift = 2pt]#1.north west) rectangle ([yshift = -6pt, xshift = -2pt]#2.south east);
  \end{tikzpicture}
}
\newcommand{\Drawb}[3]{%
  \begin{tikzpicture}[overlay,remember picture]
    \draw[#3] ([yshift = 16pt, xshift = 0pt]#1.north west) rectangle ([yshift = -1pt, xshift = 15pt]#2.south east);
  \end{tikzpicture}
}
\newcommand{\Drawc}[3]{%
  \begin{tikzpicture}[overlay,remember picture]
    \draw[#3] ([yshift = 8pt, xshift = 2pt]#1.north west) rectangle ([yshift = -3pt, xshift = 1pt]#2.south east);
  \end{tikzpicture}
}
\newcommand{\Drawd}[3]{%
  \begin{tikzpicture}[overlay,remember picture]
    \draw[#3] ([yshift = 6pt, xshift = 2pt]#1.north west) rectangle ([yshift = -1pt, xshift = -1pt]#2.south east);
  \end{tikzpicture}
}

\subsection{EPIRK4s3} 

Let us start with the fourth-order EPIRK4s3 integrator \cite{Tokman17a, Tokman17b}, the equations of which read
\begin{align}
	a^n & = u^n + \quad \tikzmark{Begin_1} \tikzmark{Begin_2} \frac{1}{8} \varphi_1\left(\frac{1}{8} \mathcal{J}(u^n)  \Delta t \right) f(u^n) \Delta t, \nonumber  \\
	b^n & = u^n + \quad \frac{1}{9} \varphi_1\left(\frac{1}{9} \mathcal{J}(u^n) \Delta t \right) f(u^n) \Delta t, \tikzmark{End_1} \nonumber  \\
	u_3^{n + 1} & = u^n + \quad \tikzmark{Begin_3} {\varphi_1\left(\mathcal{J}(u^n) \Delta t\right) f(u^n) \Delta t \tikzmark{End_2} \qquad + \tikzmark{Begin_4} {\varphi_3(\mathcal{J}(u^n) \Delta t)\left(1892 \mathcal{R}(a^n) + 1458\,(\mathcal{R}(b^n) - 2 \mathcal{R}(a^n)) \right) \Delta t,}} \tikzmark{End_3} \tikzmark{End_4} \nonumber \\
	u_4^{n + 1} & = u_3^{n + 1} + \quad \varphi_4 (\mathcal{J}(u^n) \Delta t) \left(-42336 \mathcal{R}(a^n) - 34992\,(\mathcal{R}(b^n) - 2 \mathcal{R}(a^n)) \right) \Delta t.
	\label{eq:epirk4s3}
\Drawa{Begin_1}{End_1}{red}
\Drawb{Begin_2}{End_2}{blue}
\Drawc{Begin_3}{End_3}{red}
\Drawd{Begin_4}{End_4}{blue}
\end{align}
For simplicity and brevity, we will refer to the computation of $a_n$ and $b_n$ collectively as computation of the internal stages, $u^{n+1}_3$ as the third stage, and the $\varphi_4$ function as the computation of the error estimate.

In the case of the Leja interpolation method, we compute the internal stages and $\varphi_1\left(\mathcal{J}(u^n) \Delta t \right)f(u^n)$ vertically. In order to generate the error estimate, one needs to compute $\varphi_4 (\mathcal{J}(u^n) \Delta t) \left(-42336 \mathcal{R}(a^n) - 34992\,(\mathcal{R}(b^n) - 2 \mathcal{R}(a^n)) \right) \Delta t$ separately, which necessitates that $\varphi_3(\mathcal{J}(u^n) \Delta t)\left(1892 \mathcal{R}(a^n) + 1458\,(\mathcal{R}(b^n) - 2 \mathcal{R}(a^n)) \right) \Delta t$ is also computed individually. This implementation is illustrated in blue in Eq. \eqref{eq:epirk4s3}. The internal stages, with \kiops, are computed vertically and the third stage is computed horizontally. This is depicted in red in Eq. \eqref{eq:epirk4s3}.

\begin{figure}[!tb]
    \centering
	\includegraphics[width = \textwidth]{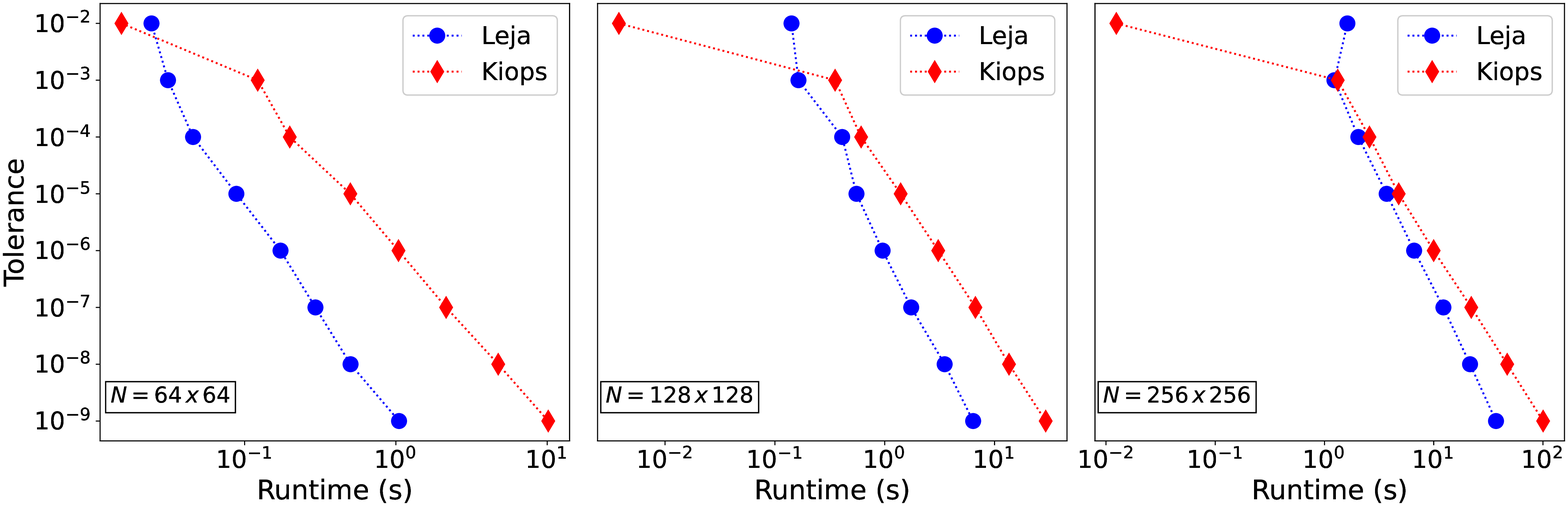} \\
	\includegraphics[width = \textwidth]{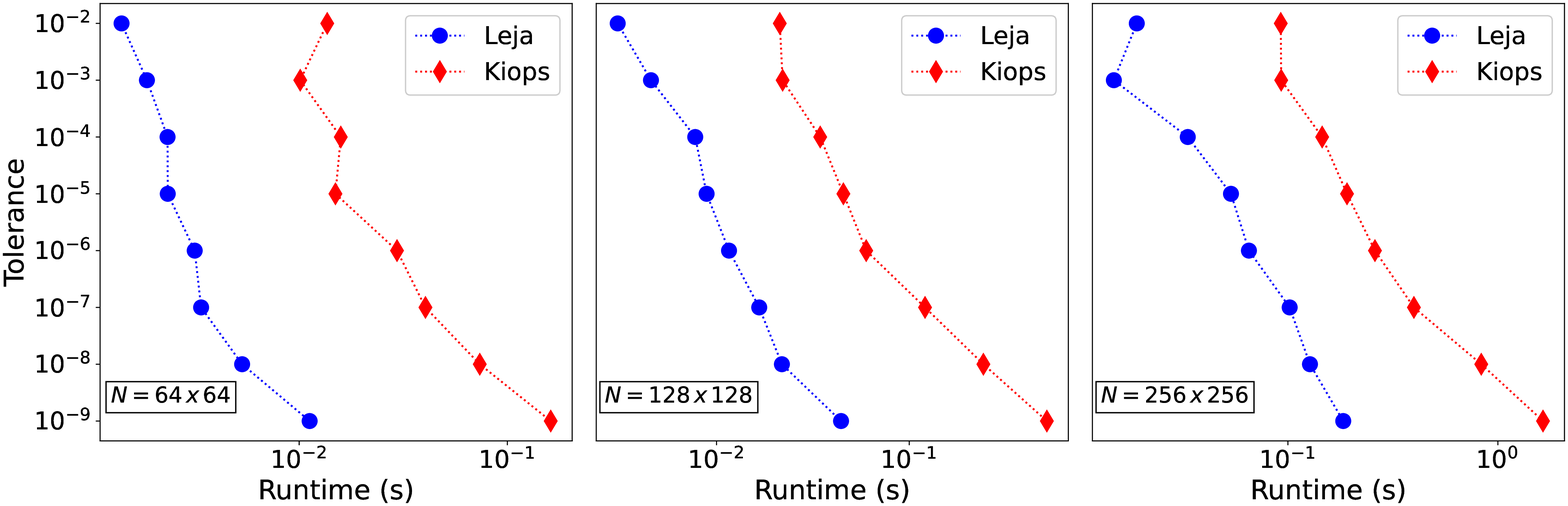}
    \caption{Comparison of the runtimes of the the functions evaluated in vertical (3 terms for Leja and 2 for \kiops) for Brusselator ($\alpha = 10^{-3}$, top panel) and Gray--Scott ($\alpha_1, \alpha_2 = 10^{-3}$,bottom panel).}
    \label{fig:leja_kiops_vertical}
\end{figure}

\begin{figure}
	\centering
	\includegraphics[width = \textwidth]{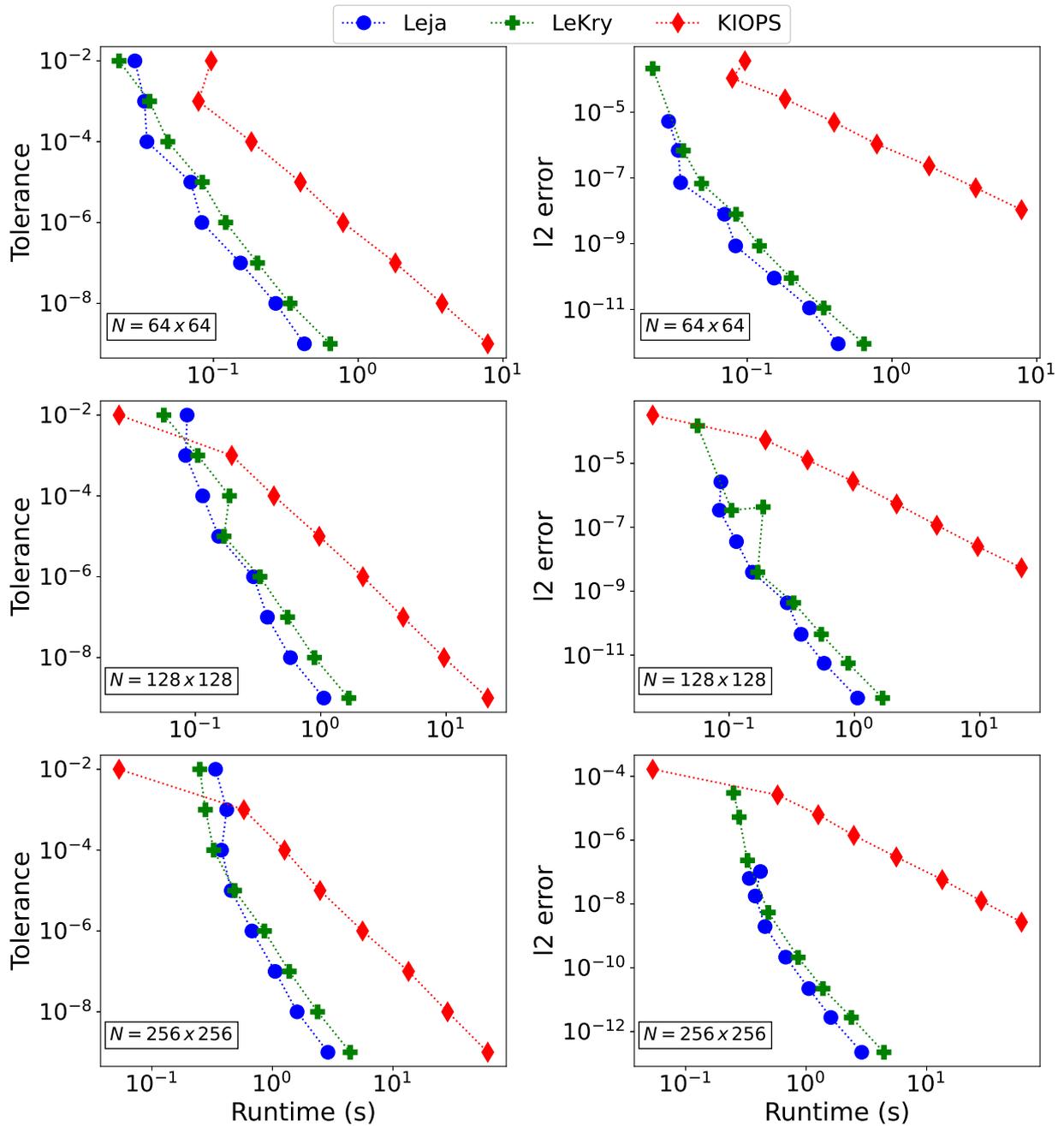}
	\caption{Comparison of the runtimes for a given value of tolerance (left panel) and the l2 norm of the global error incurred (right panel) for Leja (blue circles), LeKry (green pluses), and \kiops (red diamonds) with the EPIRK4s3 integrator for the ADR equation ($\alpha = 0.01$).}
	\label{fig:adr_4s3_0.01}
\end{figure}

\begin{figure}
	\centering
	\includegraphics[width = \textwidth]{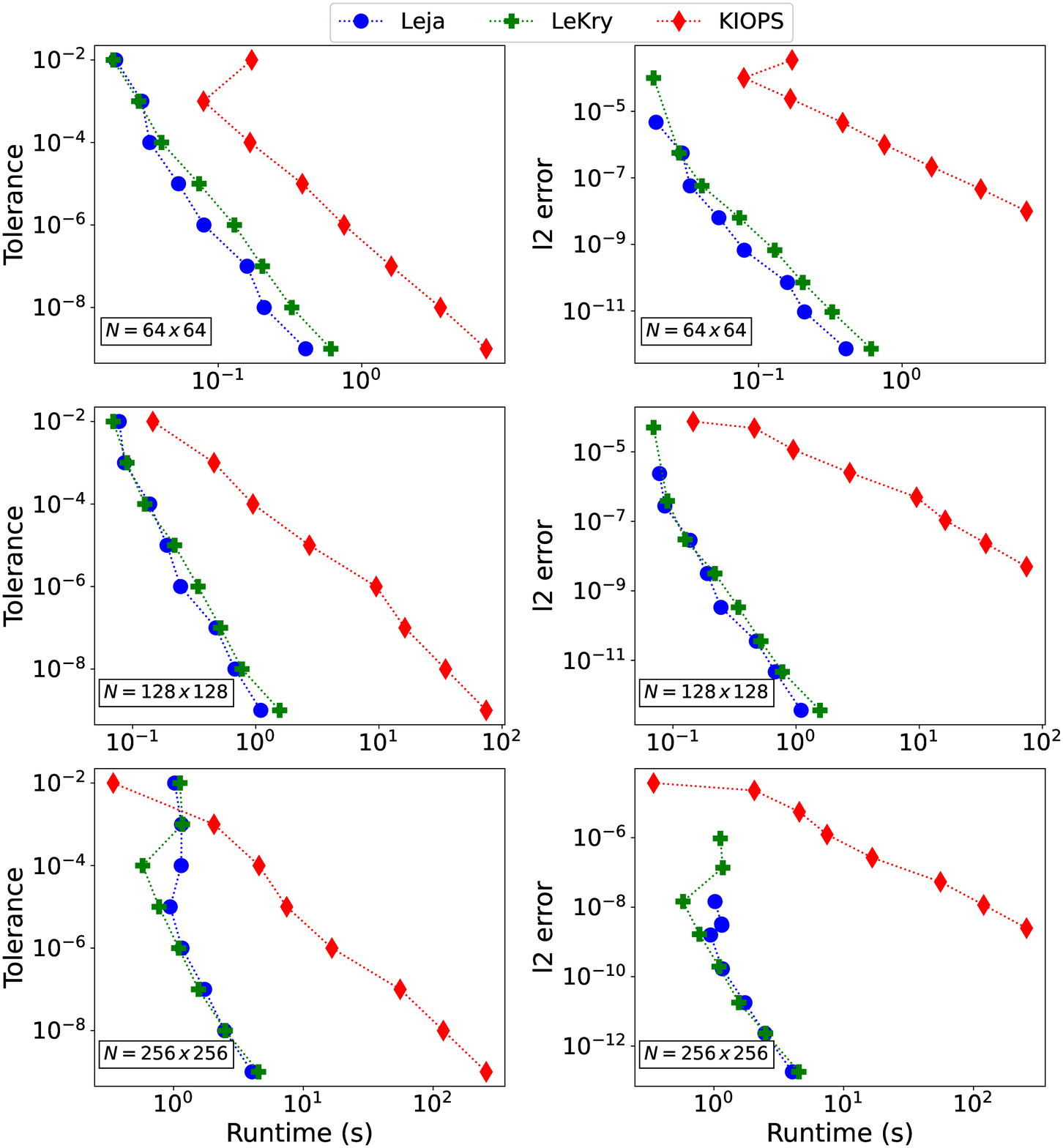} 
	\caption{Comparison of the runtimes for a given value of tolerance (left panel) and the l2 norm of the global error incurred (right panel) for Leja (blue circles), LeKry (green pluses), and \kiops (red diamonds) with the EPIRK4s3 integrator for the ADR equation ($\alpha = 0.1$).}
	\label{fig:adr_4s3_0.1}
\end{figure}

\begin{figure}
    \centering
	\includegraphics[width = \textwidth]{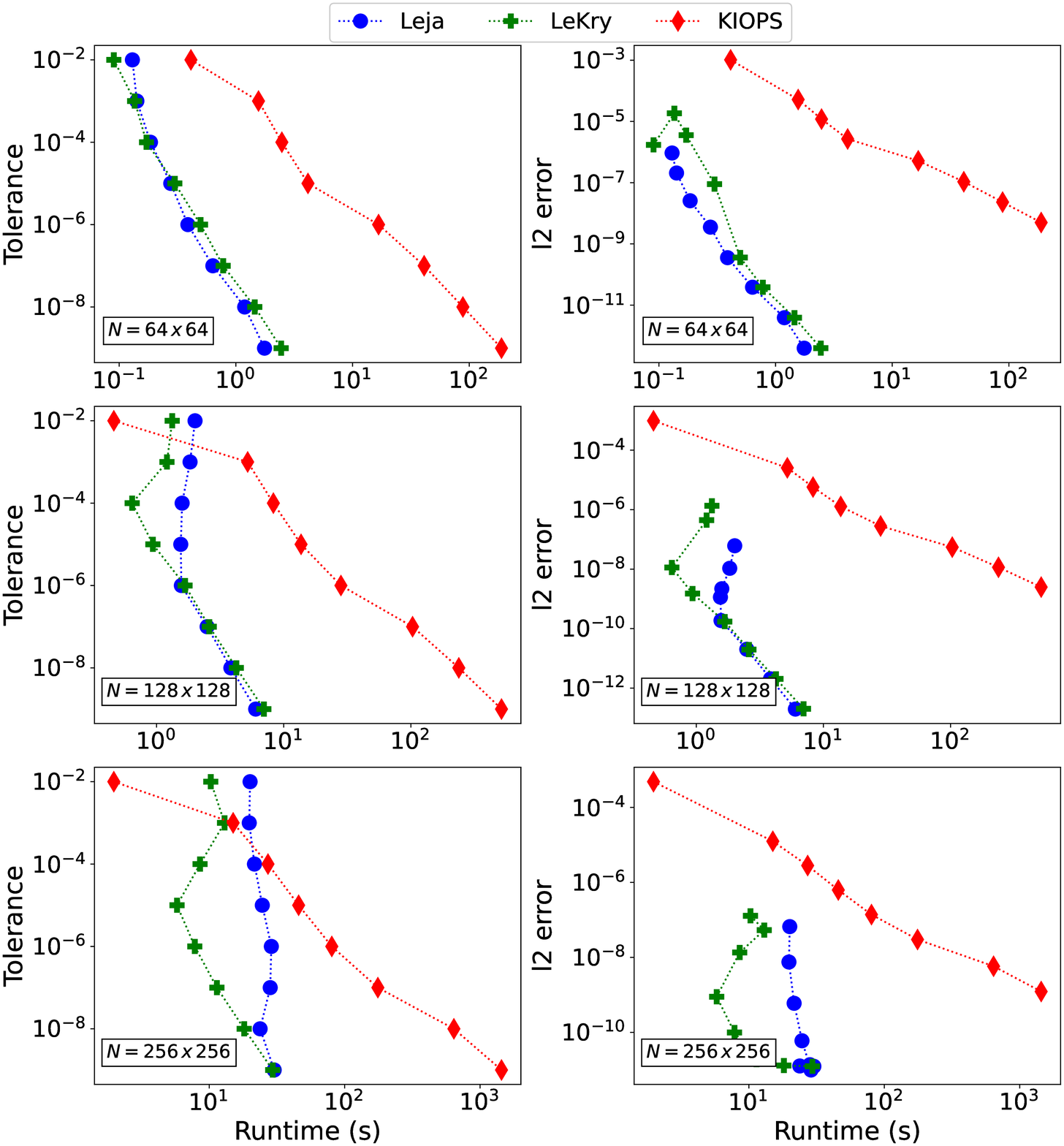} 
    \caption{Comparison of the runtimes for a given value of tolerance (left panel) and the l2 norm of the global error incurred (right panel) for Leja (blue circles), LeKry (green pluses), and \kiops (red diamonds) with the EPIRK4s3 integrator for the Brusselator equation ($\alpha = 0.001$).}
    \label{fig:bru_4s3_0.001}
\end{figure}

\begin{figure}
    \centering
	\includegraphics[width = \textwidth]{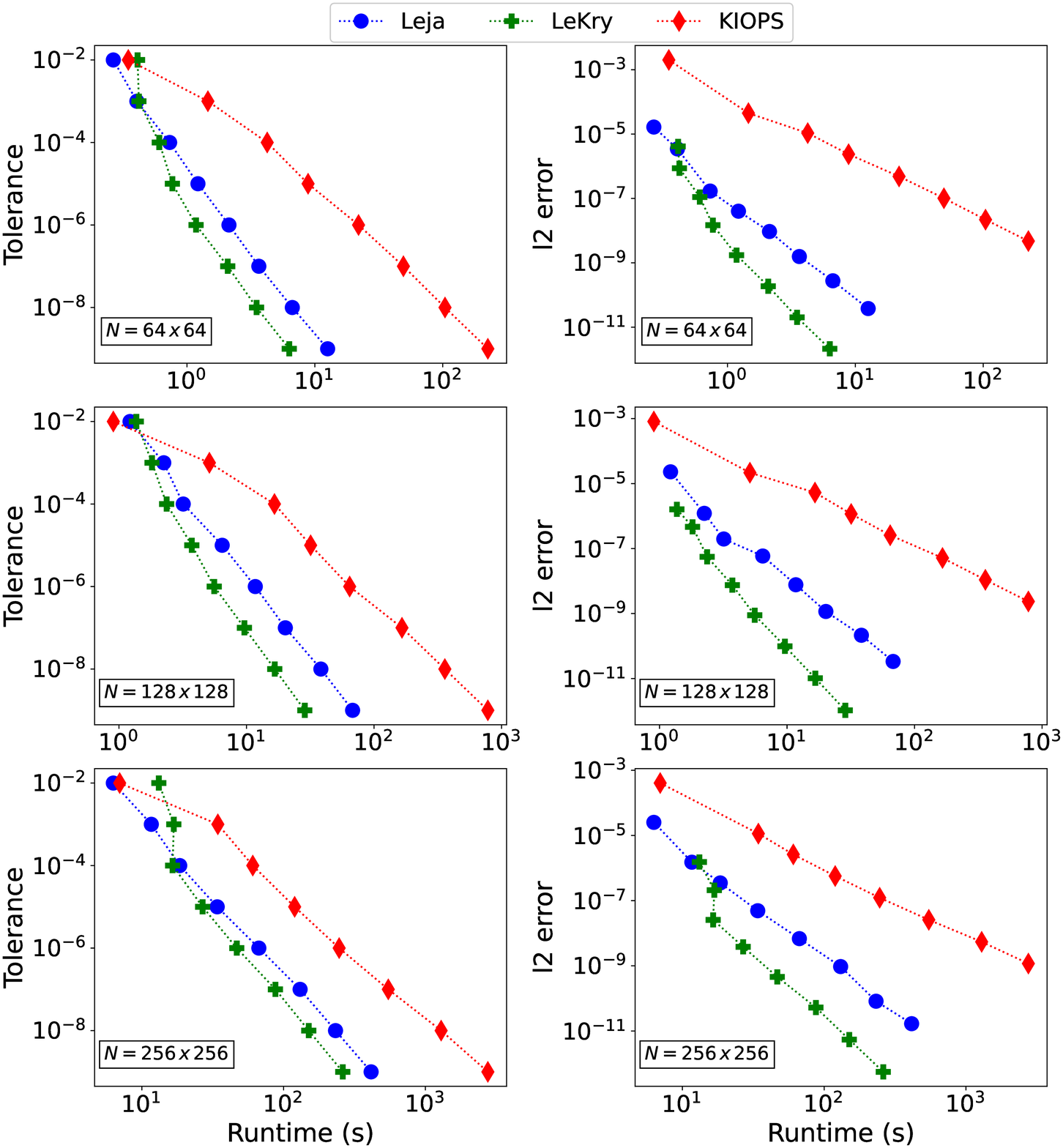} 
    \caption{Comparison of the runtimes for a given value of tolerance (left panel) and the l2 norm of the global error incurred (right panel) for Leja (blue circles), LeKry (green pluses), and \kiops (red diamonds) with the EPIRK4s3 integrator for the Brusselator equation ($\alpha = 0.1$).}
    \label{fig:bru_4s3_0.1}
\end{figure}

To understand the performance difference of the Leja method with that of \kiops in the vertical implementation, we compare, in Fig. \ref{fig:leja_kiops_vertical}, the cost of the vertical stages. The most striking feature is that even though the Leja method has to compute an additional vector, it still manages to outperform \kiops for a majority of the considered configurations. From this observation, one could extrapolate that computing just two vectors (the internal stages) with Leja, in vertical, would be even cheaper. The major advantage of this integrator over the others (considered below) is that the coefficients of the internal stages, $a_n$ and $b_n$ are rather small. This results in scaling down of the Jacobian by a substantial amount (which is analogous to choosing a smaller value of $\Delta t$), thereby expediting the convergence. This, and the fact that the horizontal implementation of \kiops, used in computing $u_3^{n+1}$, is more efficient than Leja, we propose to compute the $\varphi_l$ functions in the internal stages vertically using the Leja method whilst computing the third stage horizontally using \kiops. The error estimate, $\varphi_4 (\mathcal{J}(u^n) \Delta t) b_k$, can be computed either using the Leja or Krylov method. We compare this novel algorithm, which we call ``LeKry'' with Leja and \kiops for a wide range of problems and their corresponding parameters (Figs. \ref{fig:adr_4s3_0.01}, \ref{fig:adr_4s3_0.1}, \ref{fig:bru_4s3_0.001}, \ref{fig:bru_4s3_0.1}, and Figs. 1 - 8 in the associated supplementary material).

A comparison of these three iterative schemes, for this integrator, shows that in all considered cases, Leja and LeKry categorically perform better than \kiops in terms of computational cost as well as error incurred. In the case of the semilinear (Fig. 8 in the supplementary material) and Brusselator (with $\alpha = 0.1$, Fig. \ref{fig:bru_4s3_0.1}) equations, there is a noticeable improvement in the global error incurred by LeKry over the Leja method. In all other cases, the error incurred by both these schemes are roughly the same and is smaller than that of \kiops. It is to be noted that with the increase in the number of grid points and the magnitude of the diffusion coefficients(s), for the same user-defined tolerance, LeKry is able to better constrain the global error whilst being as expensive as, if not less than the Leja method. As such, with a reasonable degree of certainty, one can say the LeKry is the best iterative scheme, for EPIRK4s3, amongst the three.

Comparing the performance of the iterative schemes under consideration involves two factors: the computational runtime and the error incurred. The purpose of setting a tolerance for variable step size implementation is to ensure that the error, at every time step, is constrained by the tolerance. How much is the error incurred, provided that it is below the tolerance, depends on the integrator, the iterative scheme, the problem under consideration, the initial and boundary conditions, and the simulation time. As it is only the tolerance that is directly under the control of the user, we compare the runtimes of the different iterative schemes for a given tolerance. The plots for error incurred versus the runtime is used to compare the accuracy of the different iterative schemes, and not the speed-ups achieved by one over the other.

We note that comparing the different schemes at a single time step for a given step size is not a justified measure. In doing so, one ignores a plethora of other factors that include the varying spectrum of the Jacobian at every time step, the effect of the (varying) step sizes, the norm of the function being interpolated or projected, and the differences in the computation of the error incurred. For instance, say at a given time step during the simulations, the difference between the third- and fourth-order solutions of EPIRK4s3 (i.e. the error estimate) computed by the individual schemes vary at the seventh decimal place. This difference would result in the step size controller yielding different values of the step size for the next time step for the two schemes. Over a period of a few time steps, this may even translate to a different number of total time steps. As such, one needs to take into account the aforementioned factors whilst comparing the performance of the iterative schemes. Additionally, the Leja method requires an estimate of the largest eigenvalue. The time consumed by the power iterations algorithm, based on Arnoldi iterations \citep{Arnoldi1951, Vorst87}, to get the largest eigenvalue, needs to be taken into account.

\subsection{EPIRK4s3A} 

\newcommand{\Drawr}[3]{%
  \begin{tikzpicture}[overlay,remember picture]
    \draw[#3] ([yshift = 13pt, xshift = 2pt]#1.north west) rectangle ([yshift = -8pt, xshift = 1pt]#2.south east);
  \end{tikzpicture}
}
\newcommand{\Draws}[3]{%
  \begin{tikzpicture}[overlay,remember picture]
    \draw[#3] ([yshift = 11.5pt, xshift = 2pt]#1.north west) rectangle ([yshift = -6.5pt, xshift = -1pt]#2.south east);
  \end{tikzpicture}
}

\begin{align}
	a^n & = u^n + \quad \tikzmark{Begin_1} \tikzmark{Begin_2} \frac{1}{2} \varphi_1\left(\frac{1}{2} \mathcal{J}(u^n)  \Delta t \right) f(u^n) \Delta t, \nonumber  \\
	b^n & = u^n + \quad \frac{2}{3} \varphi_1\left(\frac{2}{3} \mathcal{J}(u^n) \Delta t \right) f(u^n) \Delta t, \tikzmark{End_1} \nonumber  \\
	u_3^{n + 1} & = u^n + \quad \tikzmark{Begin_3} {\varphi_1\left(\mathcal{J}(u^n) \Delta t\right) f(u^n) \Delta t \tikzmark{End_2} \qquad + \tikzmark{Begin_4}{\varphi_3(\mathcal{J}(u^n) \Delta t)\left(32 \mathcal{R}(a^n) - \frac{27}{2} \mathcal{R}(b^n)) \right) \Delta t,}} \tikzmark{End_3} \tikzmark{End_4} \nonumber \\
	u_4^{n + 1} & = u_3^{n + 1} + \quad \varphi_4 (\mathcal{J}(u^n) \Delta t) \left(-144 \mathcal{R}(a^n) + 81 \mathcal{R}(b^n) \right) \Delta t.
	\label{eq:epirk4s3a}
\Drawa{Begin_1}{End_1}{red}
\Drawb{Begin_2}{End_2}{blue}
\Drawr{Begin_3}{End_3}{red}
\Draws{Begin_4}{End_4}{blue}
\end{align}

\begin{figure}
    \centering
	\includegraphics[width = \textwidth]{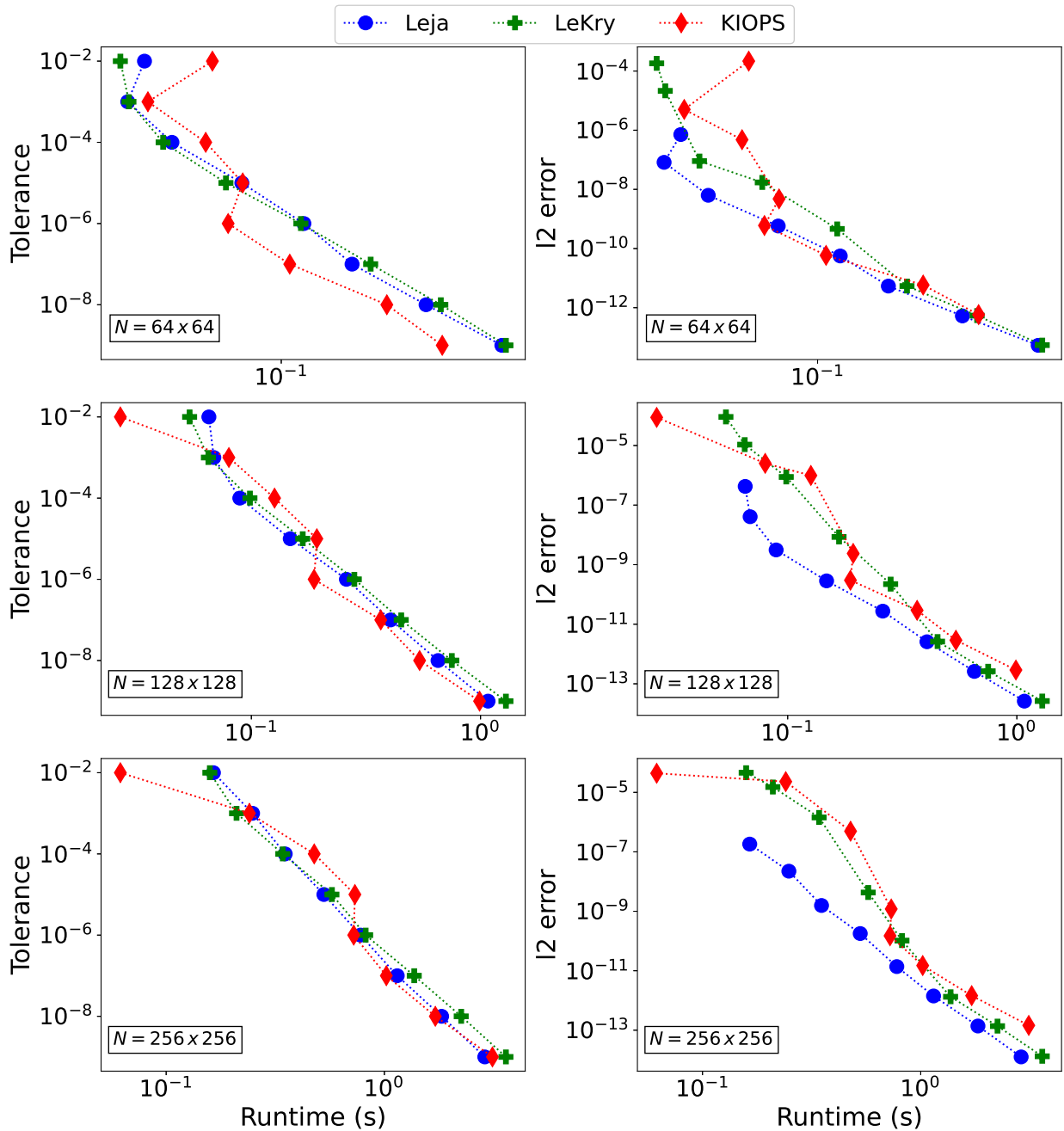} 
    \caption{Comparison of the runtimes for a given value of tolerance (left panel) and the l2 norm of the global error incurred (right panel) for Leja (blue circles), LeKry (green pluses), and \kiops (red diamonds) with the EPIRK4s3A integrator for the Allen--Cahn equation ($\alpha = 0.001$).}
    \label{fig:ac_4s3a_0.001}
\end{figure}

\begin{figure}
    \centering
	\includegraphics[width = \textwidth]{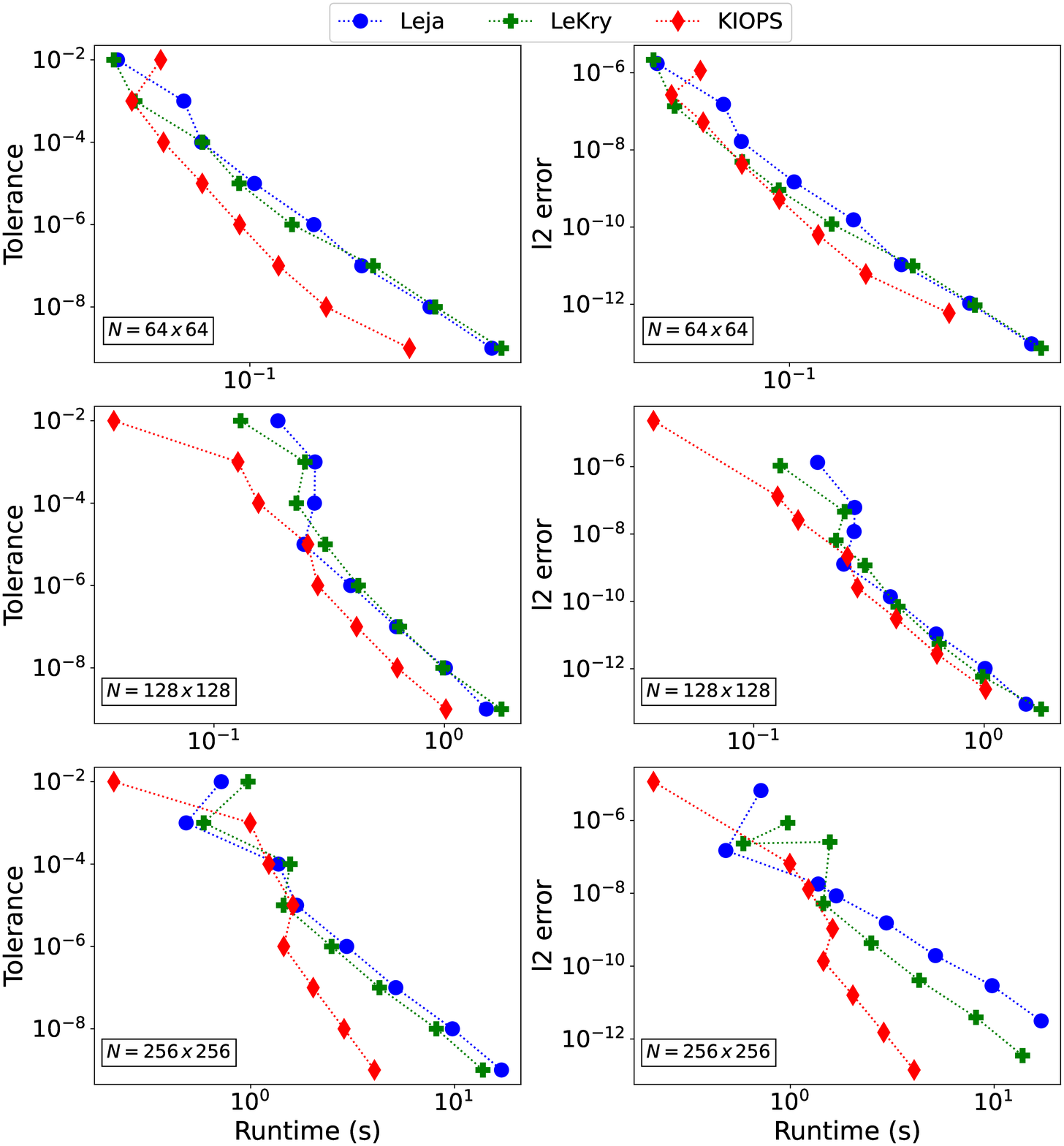} 
    \caption{Comparison of the runtimes for a given value of tolerance (left panel) and the l2 norm of the global error incurred (right panel) for Leja (blue circles), LeKry (green pluses), and \kiops (red diamonds) with the EPIRK4s3A integrator for the Allen--Cahn equation ($\alpha = 0.01$).}
    \label{fig:ac_4s3a_0.01}
\end{figure}

\begin{figure}
    \centering
	\includegraphics[width = \textwidth]{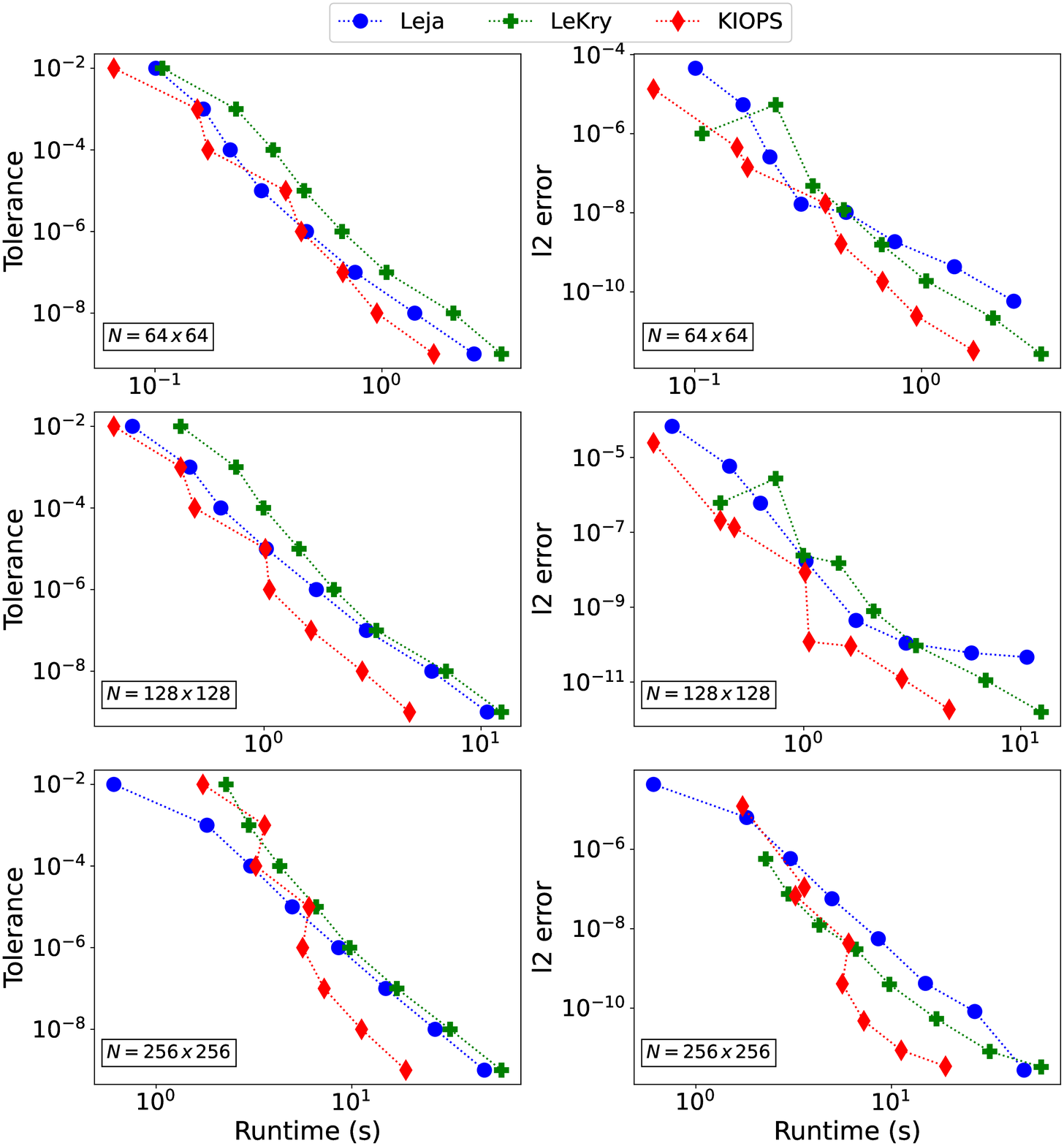} 
    \caption{Comparison of the runtimes for a given value of tolerance (left panel) and the l2 norm of the global error incurred (right panel) for Leja (blue circles), LeKry (green pluses), and \kiops (red diamonds) with the EPIRK4s3 integrator for the Allen--Cahn equation ($\alpha = 0.1$).}
    \label{fig:ac_4s3a_0.1}
\end{figure}

\begin{figure}
    \centering
	\includegraphics[width = \textwidth]{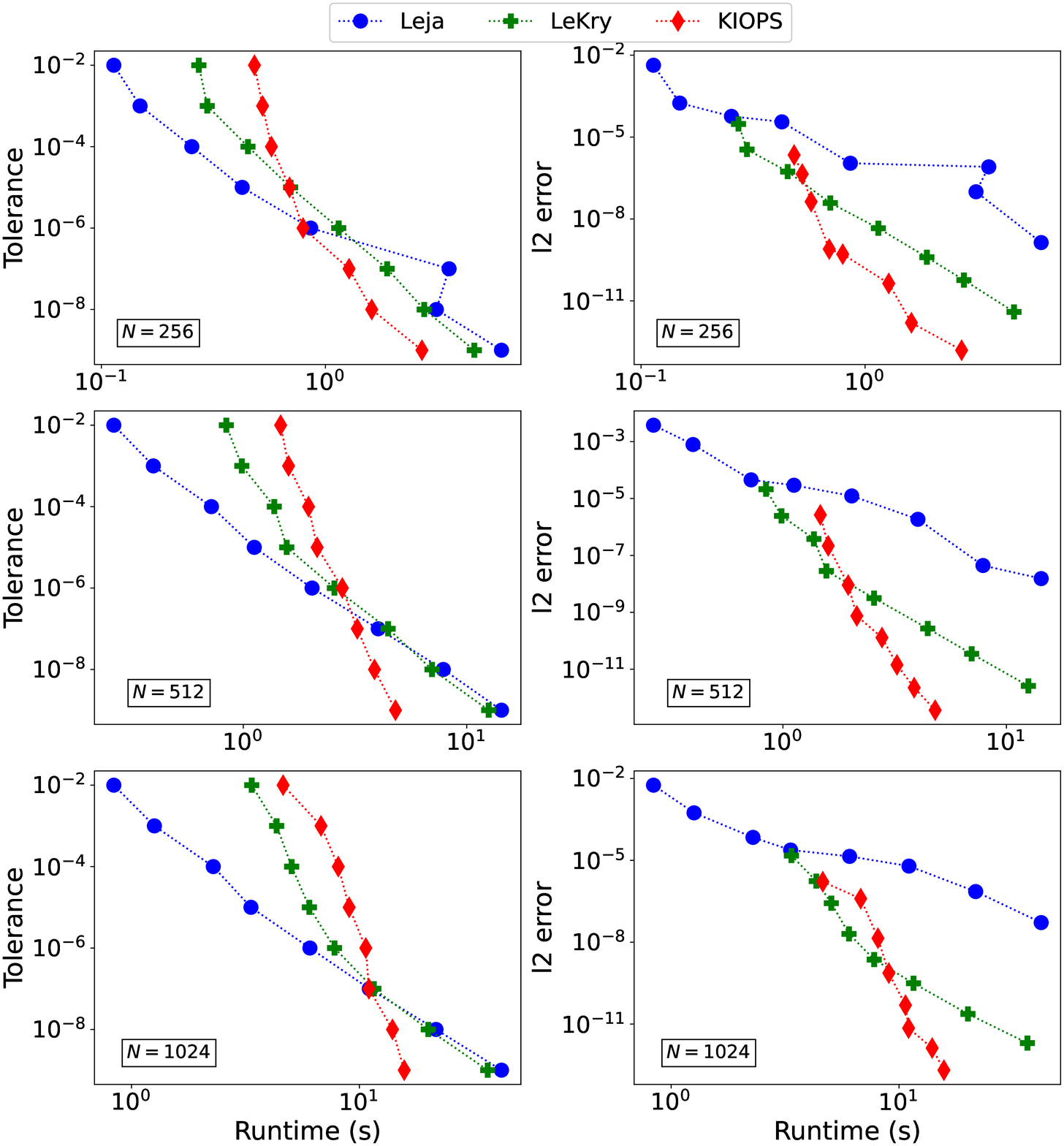} 
    \caption{Comparison of the runtimes for a given value of tolerance (left panel) and the l2 norm of the global error incurred (right panel) for Leja (blue circles), LeKry (green pluses), and \kiops (red diamonds) with the EPIRK4s3A integrator for the semilinear equation.}
    \label{fig:sl_4s3a}
\end{figure}

The EPIRK4s3A integrator (Eq. \eqref{eq:epirk4s3a}), developed by \citet{Tokman16}, is another fourth-order integrator embedded with a third-order error estimate. $a_n$ and $b_n$ are the two internal stages, $u_3^{n+1}$ is the third-order solution, and $u_4^{n+1}$ is the fourth-order solution. The order of the computation of the different $\varphi_l(z)$ functions are similar to the previous integrator: for Leja (illustrated in blue boxes), the two internal stages, and $\varphi_1\left(\mathcal{J}(u^n) \Delta t\right) f(u^n) \Delta t$ are computed vertically, and the rest individually, whereas for \kiops (shown in red boxes), the two internal stages are computed in vertical and $u_3^{n+1}$ as a linear combination of $\varphi_l(z)$ functions. As for LeKry, we compute the internal stages, in vertical, using the Leja method and $u_3^{n+1}$ using \kiops. The simulation results for a range of test problems are illustrated in Figs. \ref{fig:ac_4s3a_0.001}, \ref{fig:ac_4s3a_0.01}, \ref{fig:ac_4s3a_0.1}, \ref{fig:sl_4s3a}, and Figs. 9 - 16 in the associated supplementary material.

In the case of the Allen--Cahn equation (Figs. \ref{fig:ac_4s3a_0.001}, \ref{fig:ac_4s3a_0.01}, and \ref{fig:ac_4s3a_0.1}), Leja is somewhat faster than \kiops for the low-resolution and low-stiff cases. With the increase in the resolution and stiffness, \kiops generally performs better than Leja and LeKry. We observe similar results for the Brusselator and the Gray Scott equations (see Figs. 11 - 16 in the supplementary material). One can truly appreciate the improvements of \kiops over Leja in the case of semilinear equation (Fig. \ref{fig:sl_4s3a}) where Leja fails to converge to a reasonable accuracy for stringent tolerances. However, we see an opposite trend for the ADR equation where Leja tends to perform marginally better than LeKry and \kiops for low amounts of diffusion, whereas with the increase in the stiffness (Fig. 9 - 10 in the supplementary material), the performance of the Leja method becomes substantially better. 

We note that we do observe some differences in the computational cost for low and high tolerances for the different iterative schemes. However, as these differences are minor, we consider the overall performance of the iterative scheme for a given problem with a certain amount of stiffness. Additionally, a slight increase in the computational cost in one of the iterative schemes is usually `compensated' by an increase in the accuracy of the solution. \kiops is obviously the favoured integrator for EPIRK4s3A.

\subsection{EPIRK5P1} 

The final integrator that we consider from the class of EPIRK integrators is the fifth-order EPIRK5P1 \citep{Tokman12, Tokman13} integrator embedded with a fourth-order error estimate, the equations of which can be found in Eq. \ref{eq:epirk5p1}. The coefficients for this integrator are given in Table \ref{tab:coeff_5p1}. One particular advantage of this integrator is that all $\varphi_l(z)$ functions can be computed using the vertical implementation (for both Leja and \kiops) as follows: first, we compute $\varphi_1\left(g_{11} \, \mathcal{J}(u^n) \Delta t \right) f(u^n) \Delta t, \varphi_1\left(g_{21} \, \mathcal{J}(u^n) \Delta t \right) f(u^n) \Delta t,$ and $\varphi_1\left(g_{31} \, \mathcal{J}(u^n) \Delta t \right) f(u^n) \Delta t$ using one set of iterations. Then, we compute $\varphi_1\left(g_{22} \, \mathcal{J}(u^n) \Delta t \right) \mathcal{R}(a^n) \, \Delta t$, $\varphi_1\left(\hat{g}_{32} \, \mathcal{J}(u^n) \Delta t \right) \mathcal{R}(a^n) \, \Delta t$, and $\varphi_1\left(g_{32} \, \mathcal{J}(u^n) \Delta t \right) \mathcal{R}(a^n) \, \Delta t$ simultaneously. Once these two steps are computed, one can compute the second internal stage, i.e. $b_n$. Finally, we compute the remaining terms $\varphi_3 \left(g_{33} \, \mathcal{J}(u^n) \Delta t)(-2 \mathcal{R}(a^n) + \mathcal{R}(b^n) \right) \Delta t$ and $\varphi_3 \left(\hat{g}_{33} \, \mathcal{J}(u^n) \Delta t)(-2 \mathcal{R}(a^n) + \mathcal{R}(b^n) \right) \Delta t$. Since we do not have any linear combination of $\varphi_l(z)$ to evaluate using this procedure, we have no `LeKry' iterative scheme.

\[
	\begin{bmatrix}
	a_{11} 						\\
	a_{21} & a_{22} 			\\
	b_{1} & b_{2} &  b_{3}  	\\
	\end{bmatrix}
	=
	\begin{bmatrix}
	0.35129592695058193092 								\\
	0.84405472011657126298 		& 1.6905891609568963624 \\
	1.0 & 1.2727127317356892397 & 2.2714599265422622275 \\
	\end{bmatrix} 
\]
\[
	\begin{bmatrix}
	g_{11} 						\\
	g_{21} & g_{22} 			\\
	g_{31} & g_{32} \, (\hat{g}_{32}) & g_{33} \, (\hat{g}_{33}) 	\\
	\end{bmatrix}
	=
	\begin{bmatrix}
	0.35129592695058193092 									\\
	0.84405472011657126298 		 & 0.5 						\\
	1.0 & 0.71111095364366870359 \, (0.5) & 0.62378111953371494809 \, (1.0)	\\
	\end{bmatrix}
	\label{tab:coeff_5p1}
\]

\newcommand{\Drawf}[3]{%
  \begin{tikzpicture}[overlay,remember picture]
    \draw[#3] ([yshift = 7pt, xshift = 2pt]#1.north west) rectangle ([yshift = -2pt, xshift = 1.5pt]#2.south east);
  \end{tikzpicture}
}

\begin{align}
	a^n = u^n & + \tikzmark{Begin_1} a_{11} \, \varphi_1\left(g_{11} \, \mathcal{J}(u^n) \Delta t \right) f(u^n) \Delta t \nonumber \\
	b^n = u^n & + a_{21} \, \varphi_1\left(g_{21} \, \mathcal{J}(u^n) \Delta t \right)f(u^n) \Delta t + \tikzmark{Begin_2} a_{22} \, \varphi_1\left(g_{22} \, \mathcal{J}(u^n) \Delta t \right) \mathcal{R}(a^n) \, \Delta t \nonumber \\
	u_4^{n + 1} = u^n & + b_1 \, \varphi_1\left(g_{31} \, \mathcal{J}(u^n) \Delta t\right) f(u^n) \Delta t \tikzmark{End_1} \, + b_2 \, \varphi_1\left(\hat{g}_{32} \, \mathcal{J}(u^n) \Delta t\right) \mathcal{R}(a^n) \, \Delta t \, + \tikzmark{Begin_3} b_3 \, \varphi_3 \left(\hat{g}_{33} \, \mathcal{J}(u^n) \Delta t)(-2 \mathcal{R}(a^n) + \mathcal{R}(b^n) \right) \Delta t \nonumber \\
	u_5^{n + 1} = u^n & + b_1 \, \varphi_1\left(g_{31} \, \mathcal{J}(u^n) \Delta t\right) f(u^n) \Delta t \, +  b_2 \, \varphi_1\left(g_{32} \, \mathcal{J}(u^n) \Delta t\right) \mathcal{R}(a^n) \, \Delta t \tikzmark{End_2} \, + b_3 \, \varphi_3 \left(g_{33} \, \mathcal{J}(u^n) \Delta t)(-2 \mathcal{R}(a^n) + 	\mathcal{R}(b^n) \right) \Delta t \tikzmark{End_3}
	\Drawf{Begin_1}{End_1}{blue}
	\Drawf{Begin_2}{End_2}{blue}
	\Drawf{Begin_3}{End_3}{blue}
	\label{eq:epirk5p1}
\end{align}

\begin{figure}
    \centering
	\includegraphics[width = \textwidth]{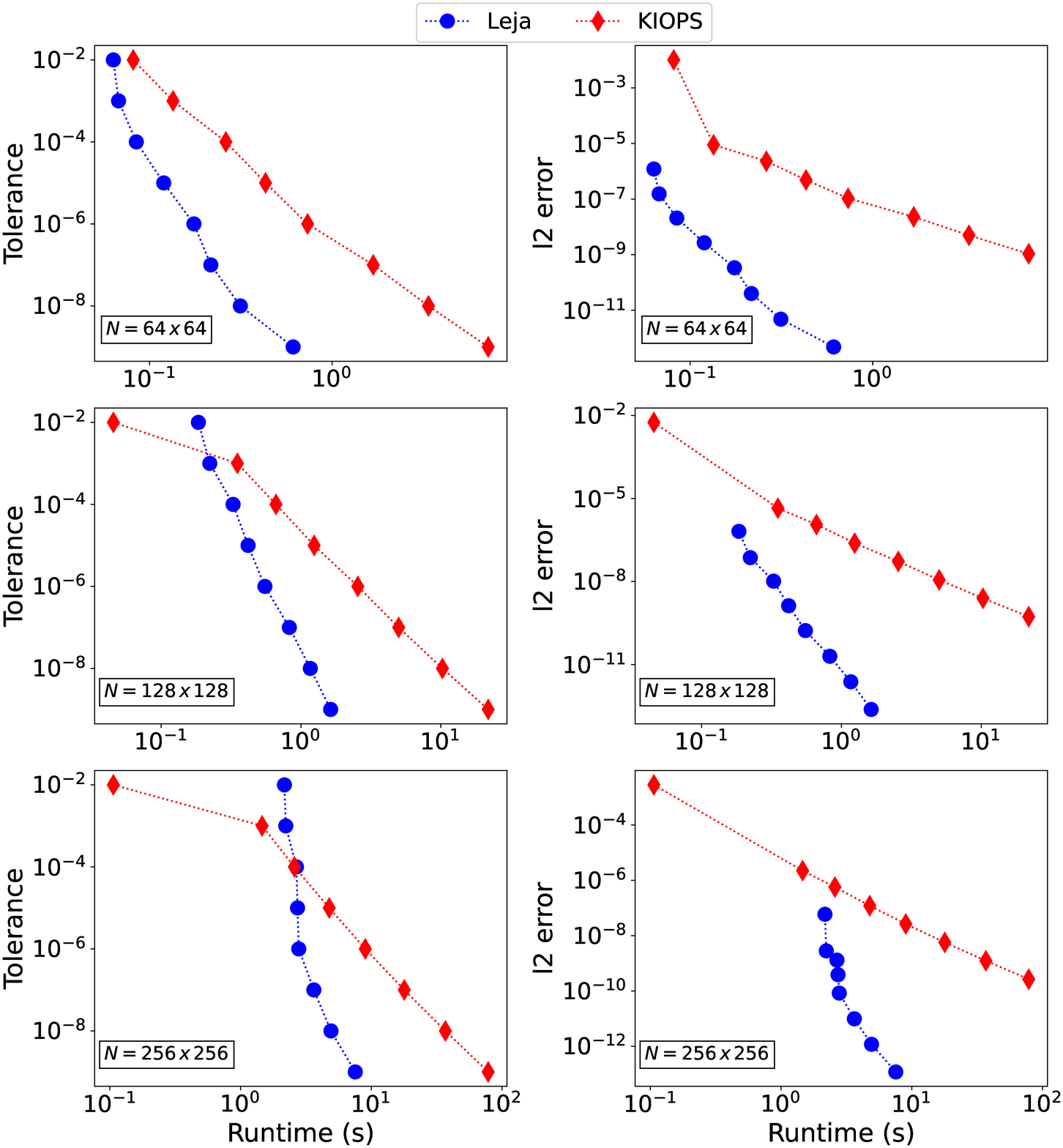} 
    \caption{Comparison of the runtimes for a given value of tolerance (left panel) and the l2 norm of the global error incurred (right panel) for Leja (blue circles), LeKry (green pluses), and \kiops (red diamonds) with the EPIRK5P1 integrator for the Brusselator equation ($\alpha = 0.001$).}
    \label{fig:bru_5p1_0.001}
\end{figure}

\begin{figure}
    \centering
	\includegraphics[width = \textwidth]{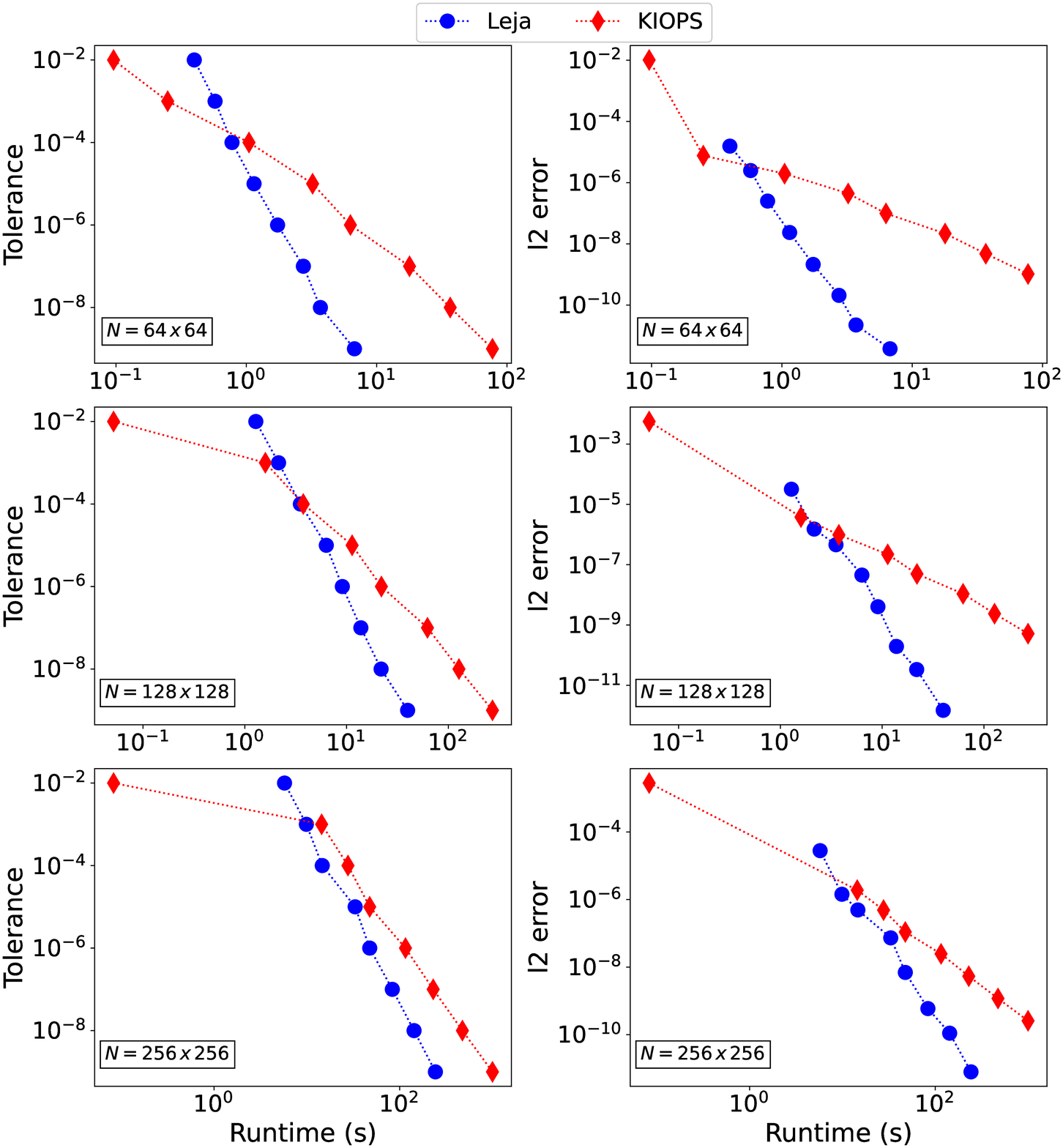} 
	\caption{Comparison of the runtimes for a given value of tolerance (left panel) and the l2 norm of the global error incurred (right panel) for Leja (blue circles), LeKry (green pluses), and \kiops (red diamonds) with the EPIRK5P1 integrator for the Brusselator equation ($\alpha = 0.1$).}
    \label{fig:bru_5p1_0.1}
\end{figure}

\begin{figure}
    \centering
	\includegraphics[width = \textwidth]{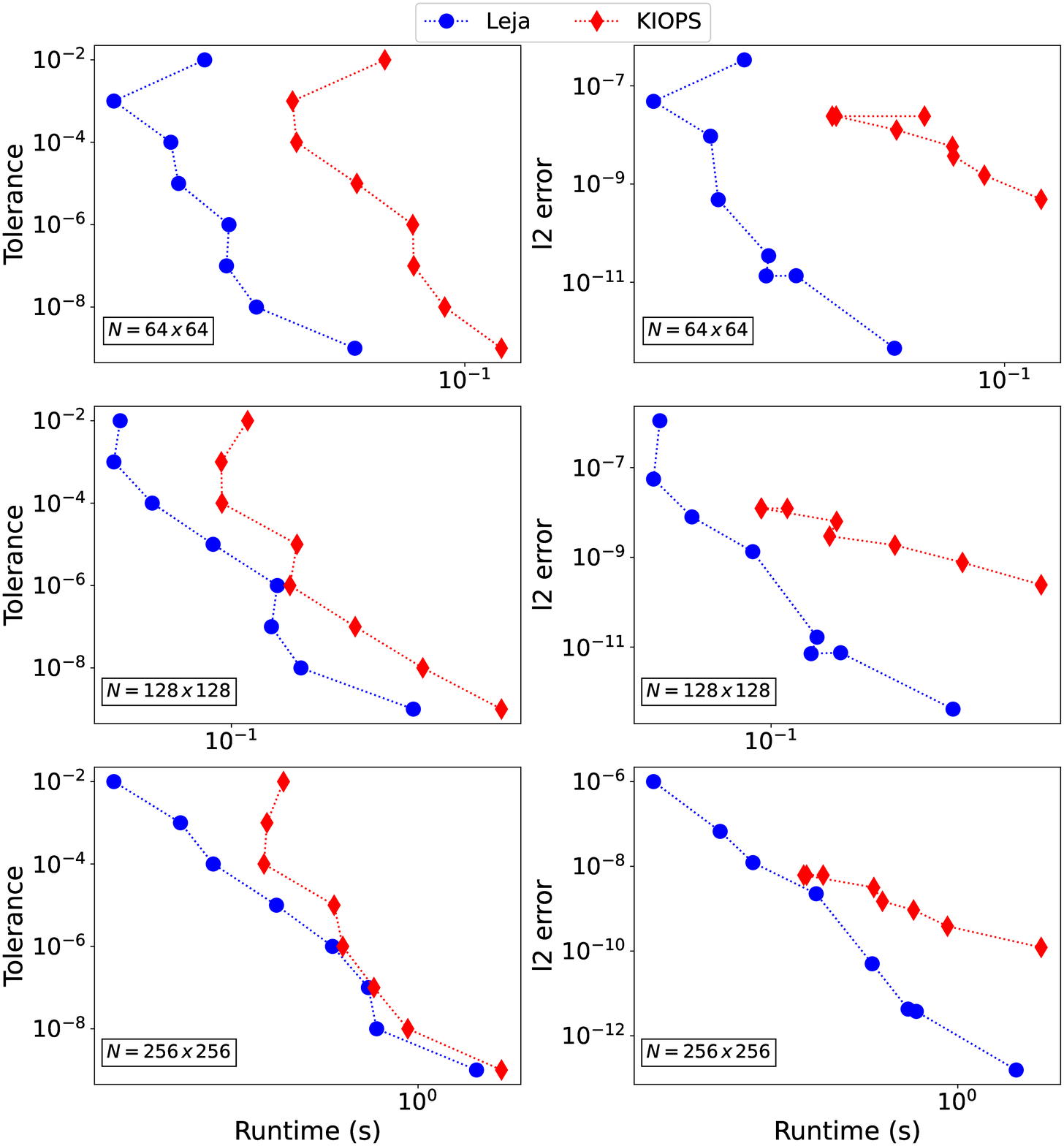} 
    \caption{Comparison of the runtimes for a given value of tolerance (left panel) and the l2 norm of the global error incurred (right panel) for Leja (blue circles), LeKry (green pluses), and \kiops (red diamonds) with the EPIRK5P1 integrator for the Gray--Scott equation ($\alpha = 0.001$).}
    \label{fig:gs_5p1_0.001}
\end{figure}

\begin{figure}
    \centering
	\includegraphics[width = \textwidth]{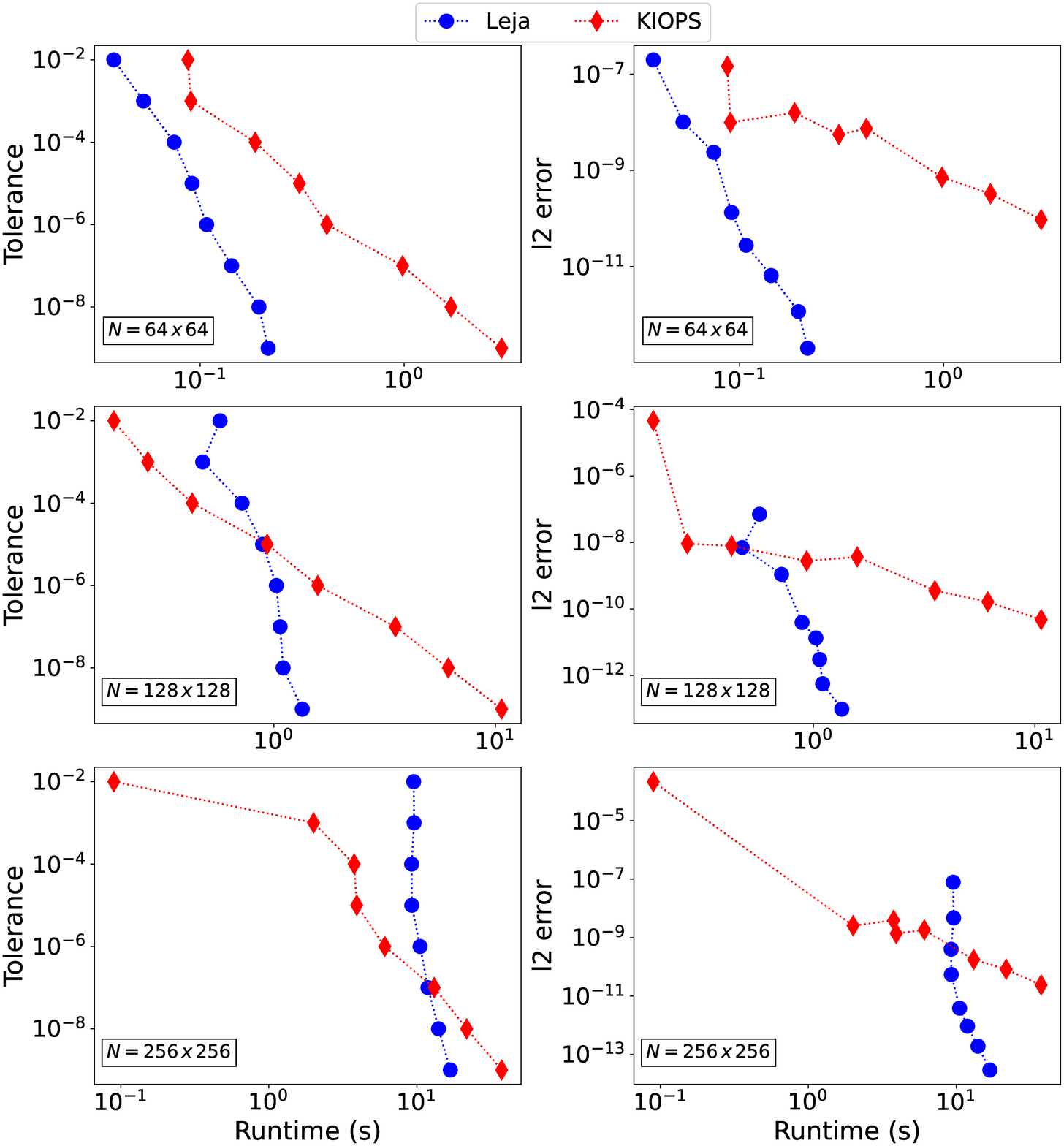} 
    \caption{Comparison of the runtimes for a given value of tolerance (left panel) and the l2 norm of the global error incurred (right panel) for Leja (blue circles), LeKry (green pluses), and \kiops (red diamonds) with the EPIRK5P1 integrator for the Gray--Scott equation ($\alpha = 0.1$).}
    \label{fig:gs_5p1_0.1}
\end{figure}

The simulation results for the set of problems using EPIRK5P1 are shown in Figs. \ref{fig:bru_5p1_0.001}, \ref{fig:bru_5p1_0.1}, \ref{fig:gs_5p1_0.001}, \ref{fig:gs_5p1_0.1}, and Figs. 17 - 24 in the supplementary material. In the case of the Brusselator equation (Figs. \ref{fig:bru_5p1_0.001} and \ref{fig:bru_5p1_0.1}), we find that Leja has superior performance than that of \kiops in almost all considered configurations. Whilst in some cases Leja has a slightly increased cost, the error incurred in those cases is consistently lower than that of \kiops. For a wide range of tolerance, Leja incurs low computational cost whilst yielding more accurate solutions for the ADR equation (supplementary material Figs. 17 - 18). With the increase in the magnitude of the diffusion coefficient and the resolution, Leja tends to be slightly more expensive than \kiops in the low tolerance regime. We see the opposite trend in the Allen--Cahn equation (supplementary material Figs. 19 - 21). For low amounts of diffusion and low resolution, \kiops has a faster runtime but Leja tends to yield more accurate solutions. As the stiffness (and resolution) is increased, the difference in the computational runtimes becomes minuscule. The results for the Gray--Scott equation (Figs. \ref{fig:gs_5p1_0.001} and \ref{fig:gs_5p1_0.1}) are similar to that of the ADR equation. For the lower amounts of diffusion, Leja has improved runtimes (and accuracy) over \kiops. However, as the amount of diffusion (and the resolution) is increased, \kiops tends to be faster than Leja in the lenient to intermediate tolerance regime. For stringent tolerances, Leja has somewhat improved runtimes whilst yielding significantly more accurate results. Finally, for the semilinear equation (supplementary material Fig. 24), we find that the difference in runtimes between the two iterative schemes is marginal; \kiops being faster in the stringent tolerance range. Overall, we find that the performance of Leja and \kiops, for EPIRK5P1, depends on the problem under consideration, the stiffness, and the resolution. Nevertheless, Leja seems to be the preferred scheme by preponderance of the evidence.

In \citet{Deka22a}, we reported better computational performance for the Krylov-based method, whilst maintaining high accuracy, for EPIRK5P1 in the case of magnetic reconnection. However, it is to be noted that in that work, we used a standard Krylov iterative scheme (not the highly advanced \kiops method). We did not have the vertical implementation (for either Leja or Krylov) in the aforementioned work, i.e., we computed all relevant $\varphi_l(z)$ functions individually. Furthermore, is likely that for the specific problem considered in that study, i.e. magnetic reconnection (and for the specific resolution and parameters), Krylov has better computational performance than that of Leja. This dependence of performance of an integrator or an iterative scheme on the problem and its parameters under study is in agreement with what we have found in this work. One or more of these reasons is likely why do not get results consistent with our previous work, with EPIRK5P1, for the set of MHD equations.


\subsection{EXPRB43} 

Now, we move on to the set of EXPRB integrators: we start off with the fourth-order EXPRB43 (embedded with a third-order error estimate, Eq. \ref{eq:exprb43}) integrator \citep{Caliari09, Ostermann10}. The computation of $\varphi_1\left(\frac{1}{2} \mathcal{J}(u^n)  \Delta t \right) f(u^n) \Delta t$ and $\varphi_1\left(\mathcal{J}(u^n)  \Delta t \right) f(u^n) \Delta t$ is performed in vertical. The remainder of the terms are computed individually for both Leja and \kiops. Similar to EPIRK5P1, we do not have a `LeKry' scheme. 
\begin{align}
	a^n & = u^n + \frac{1}{2} \varphi_1\left(\frac{1}{2} \mathcal{J}(u^n)  \Delta t \right) f(u^n) \Delta t \nonumber  \\
	b^n & = u^n + \varphi_1\left(\mathcal{J}(u^n) \Delta t \right) f(u^n) \Delta t + \varphi_1\left(\mathcal{J}(u^n)  \Delta t \right) \mathcal{R}(a^n) \Delta t \nonumber \\
	u_3^{n + 1} & = u^n + \varphi_1\left(\mathcal{J}(u^n) \Delta t\right) f(u^n) \Delta t + \varphi_3(\mathcal{J}(u^n) \Delta t)(16 \mathcal{R}(a^n) - 2 \mathcal{R}(b^n)) \Delta t \nonumber \\
	u_4^{n + 1} & = u_3^{n + 1} + \varphi_4(\mathcal{J}(u^n) \Delta t)(- 48 \mathcal{R}(a^n) + 12 \mathcal{R}(b^n)) \Delta t
	\label{eq:exprb43}
\end{align}

\begin{figure}
    \centering
	\includegraphics[width = \textwidth]{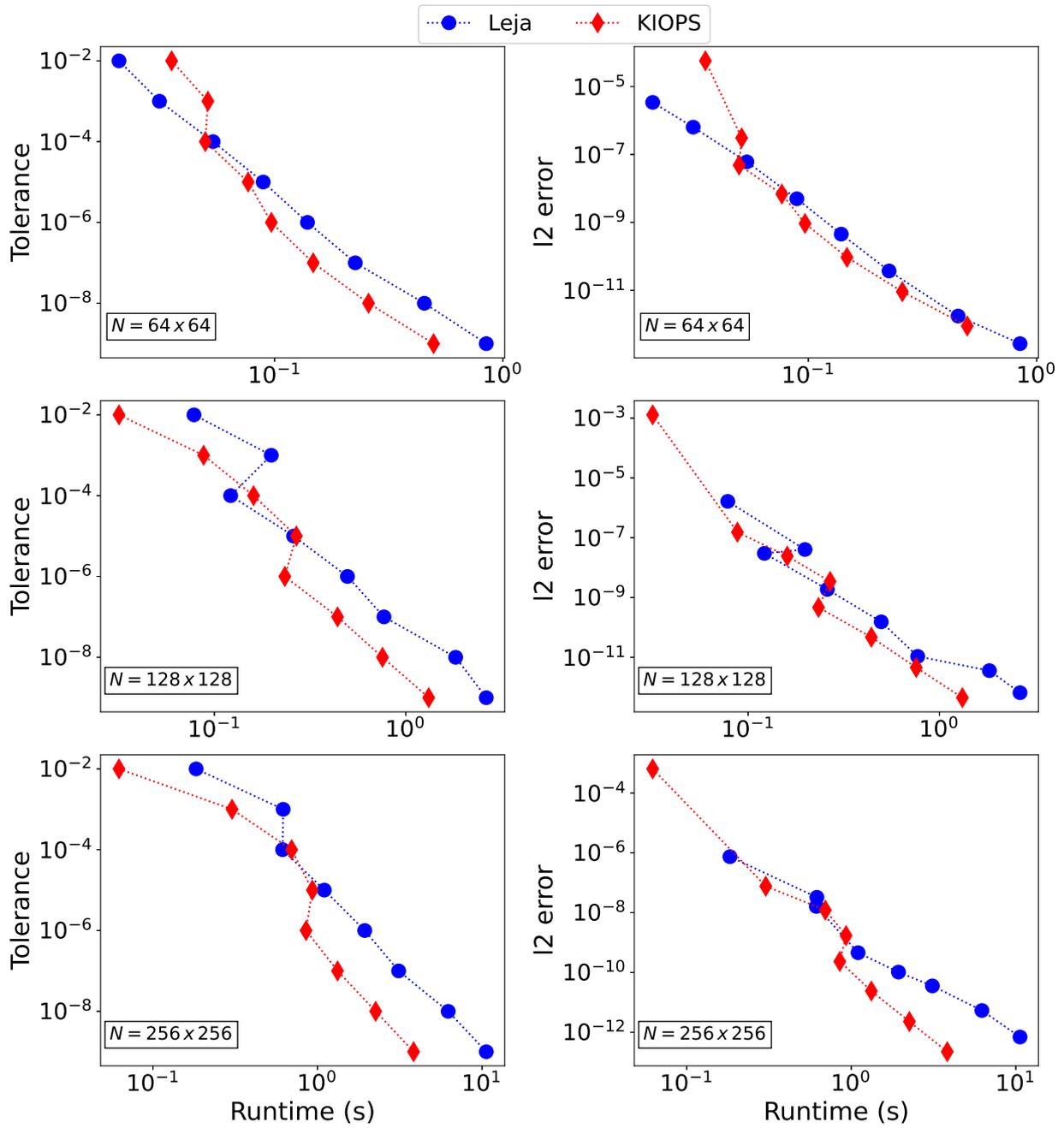} 
    \caption{Comparison of the runtimes for a given value of tolerance (left panel) and the l2 norm of the global error incurred (right panel) for Leja (blue circles), LeKry (green pluses), and \kiops (red diamonds) with the EXPRB43 integrator for the Allen--Cahn equation ($\alpha = 0.001$).}
    \label{fig:ac_43_0.001}
\end{figure}

The simulation results for the set of problems with the EXPRB43 integrator are illustrated in Figs. \ref{fig:ac_43_0.001}, \ref{fig:ac_43_0.1}, \ref{fig:bru_43_0.001}, \ref{fig:bru_43_0.1}, and Figs. 25 - 32 in the supplementary material. \kiops has similar performance to that of Leja for the Allen--Chan equation with lower amounts of diffusion (Fig. \ref{fig:ac_43_0.001}). As the stiffness in the problem is increased, \kiops clearly outperforms Leja (Figs. \ref{fig:ac_43_0.1}). Similar results can be seen for the Brusselator and the Gray--Scott equation (Figs. \ref{fig:bru_43_0.001}, \ref{fig:bru_43_0.1}, and Figs. 28 - 31 in the supplementary material). Whilst Leja might be favourable for the low stiffness and low resolution cases, \kiops is clear better of the two for stiff and high-resolution scenarios. Leja and \kiops have roughly similar runtimes for the ADR equations (supplementary material Figs. 25 - 26) with Leja yielding slightly more accurate solutions.  Finally, Leja fails to converge in the intermediate to stringent tolerance regime for the semilinear equation whereas \kiops provides highly accurate solutions with a reasonable amount of computational cost (supplementary material Figs. 32). Even though there are a few instances where Leja performs marginally better than \kiops, it is evident that \kiops is the superior iterative scheme for EXPRB43.

The conclusions that can be drawn from the results shown here are not completely in agreement with the results for the Kelvin-Helmholtz instability in \citet{Deka22a}. In that work, we find Leja to be superior to Krylov for the EXPRB43 integrator. Again, we emphasise that the improved algorithm for both Leja (vertical implementation) and Krylov (\kiops) has a significant impact on the overall performance of an integrator. Furthermore, the problem (and its associated parameters) under consideration plays a major role in the overall computational results (for example, see Figs. \ref{fig:bru_43_0.001} and supplementary material Figs. 29 and 31). 

\begin{figure}
    \centering
	\includegraphics[width = \textwidth]{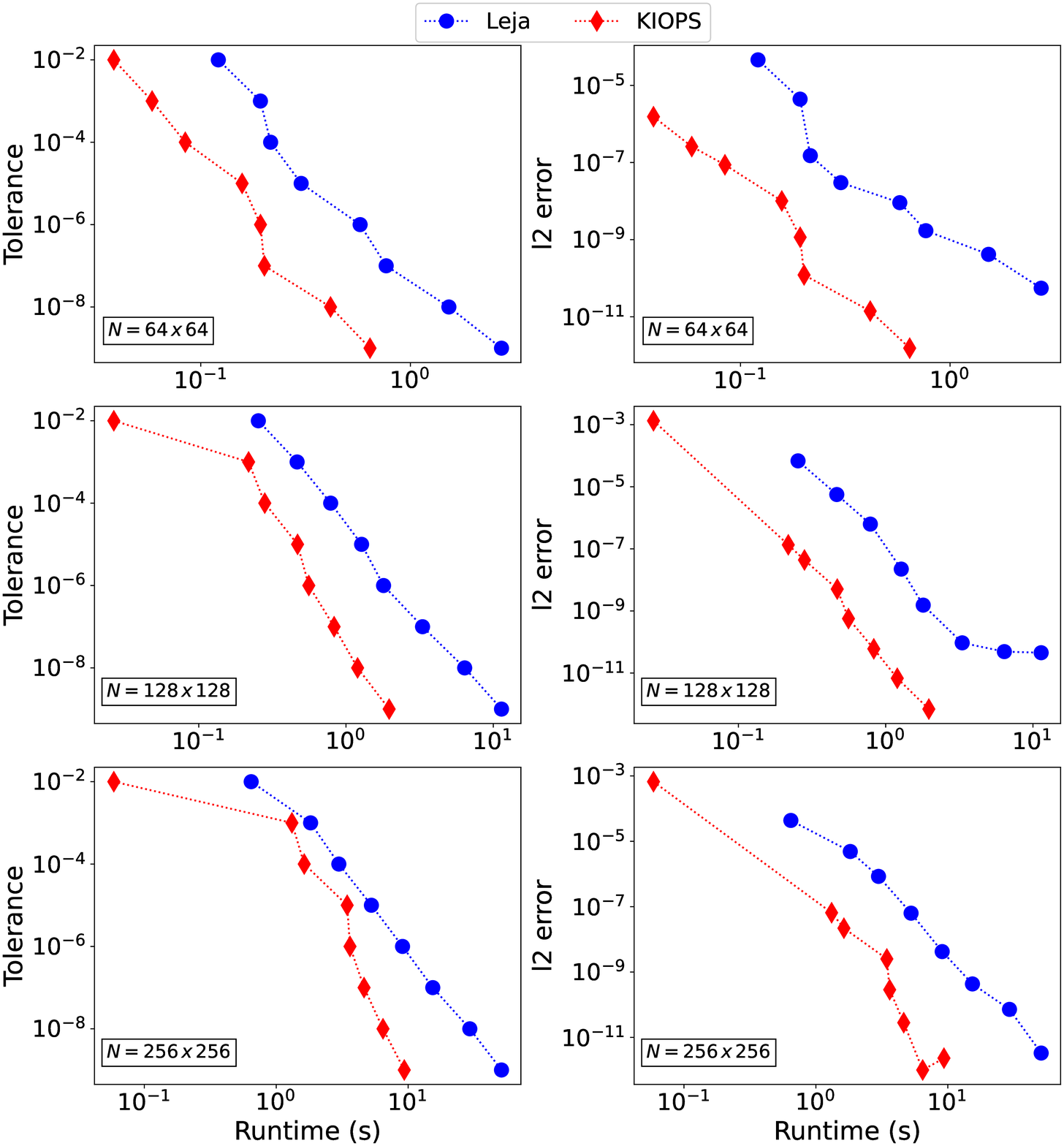} 
    \caption{Comparison of the runtimes for a given value of tolerance (left panel) and the l2 norm of the global error incurred (right panel) for Leja (blue circles), LeKry (green pluses), and \kiops (red diamonds) with the EXPRB43 integrator for the Allen--Cahn equation ($\alpha = 0.1$).}
    \label{fig:ac_43_0.1}
\end{figure}

\begin{figure}
    \centering
	\includegraphics[width = \textwidth]{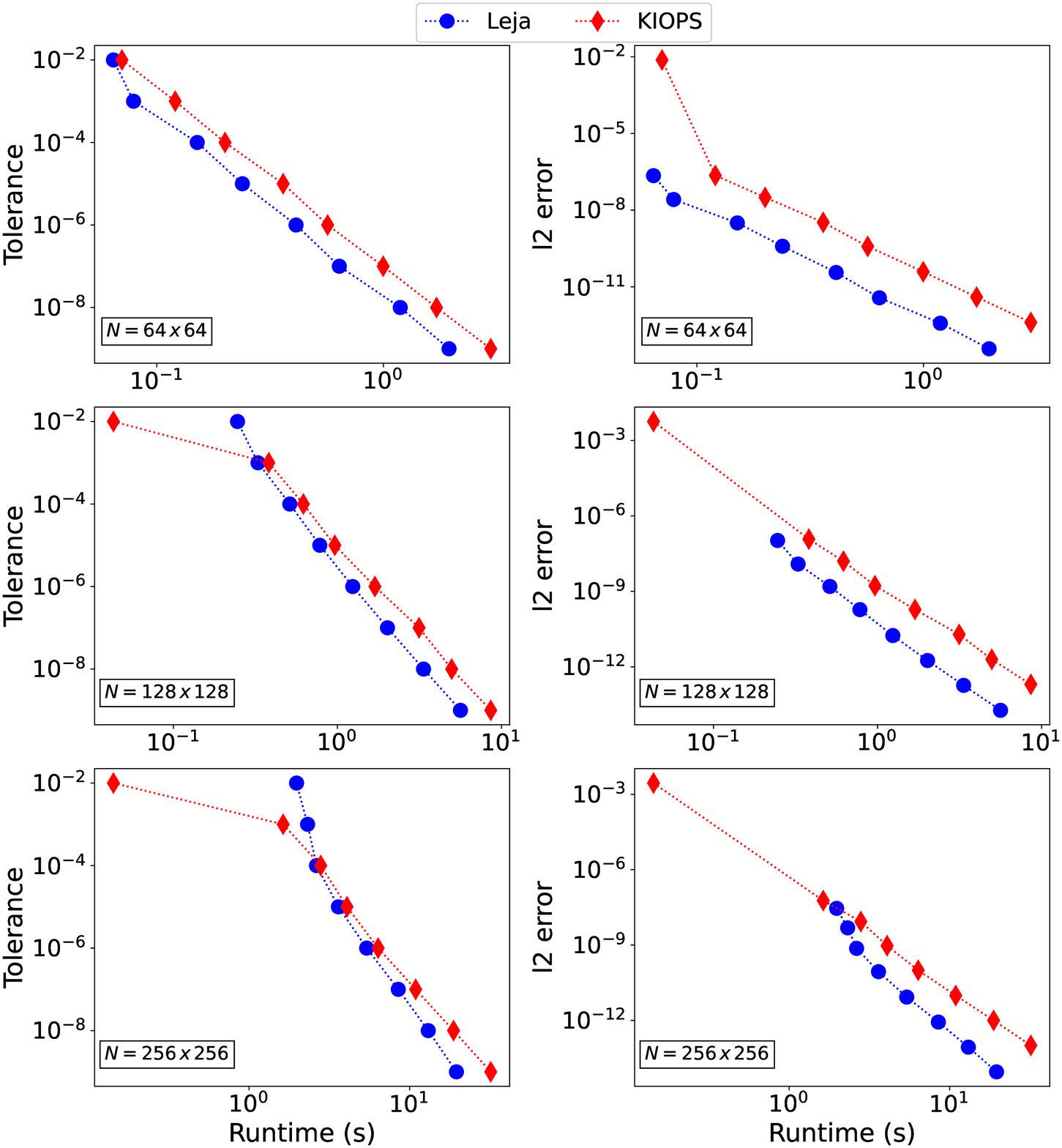} 
    \caption{Comparison of the runtimes for a given value of tolerance (left panel) and the l2 norm of the global error incurred (right panel) for Leja (blue circles), LeKry (green pluses), and \kiops (red diamonds) with the EXPRB43 integrator for the Brusselator equation ($\alpha = 0.001$).}
    \label{fig:bru_43_0.001}
\end{figure}

\begin{figure}
    \centering
	\includegraphics[width = \textwidth]{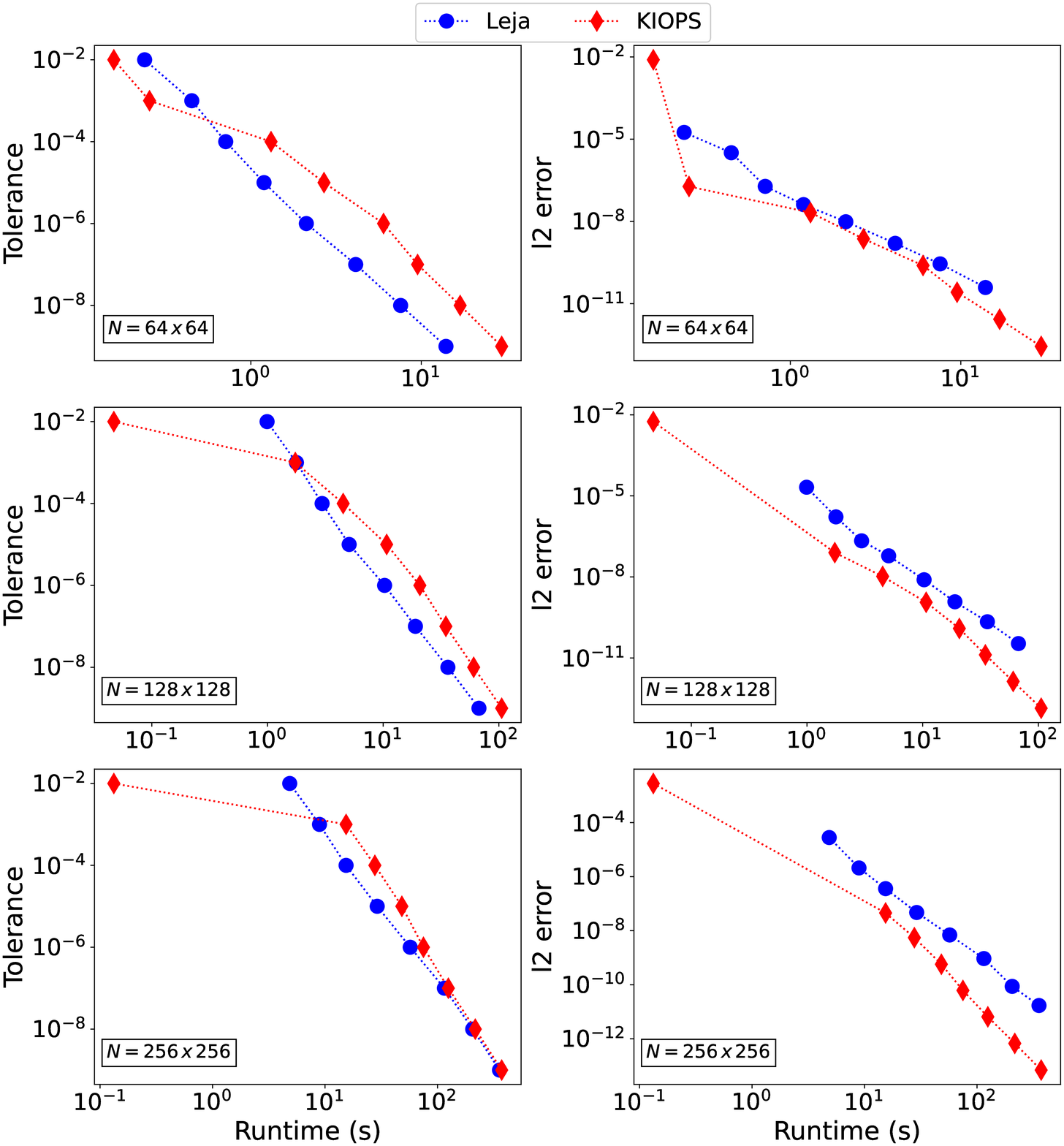}
    \caption{Comparison of the runtimes for a given value of tolerance (left panel) and the l2 norm of the global error incurred (right panel) for Leja (blue circles), LeKry (green pluses), and \kiops (red diamonds) with the EXPRB43 integrator for the Brusselator equation ($\alpha = 0.1$).}
    \label{fig:bru_43_0.1}
\end{figure}


\subsection{EXPRB53s3} 

The second EXPRB (and in this study, the last) integrator that we consider is the fifth-order EXPRB53s3 \citep{Luan14} embedded with a third-order error estimate (Eq. \eqref{eq:exprb5s3}). \kiops computes $\varphi_1\left(\frac{1}{2} \mathcal{J}(u^n) \Delta t \right) f(u^n) \Delta t$ and $\varphi_1\left(\frac{9}{10} \mathcal{J}(u^n) \Delta t \right) f(u^n) \Delta t$ as well as $\varphi_3\left(\frac{1}{2}\mathcal{J}(u^n) \Delta t\right) \mathcal{R}(a^n) \Delta t$ and $\varphi_3\left(\frac{9}{10}\mathcal{J}(u^n) \Delta t\right) \mathcal{R}(a^n) \Delta t$ in vertical. $u_3^{n + 1}$ and $u_5^{n + 1}$ are computed as linear combinations of $\varphi_l(z)$ functions, leading to an overall of four $\varphi_l(z)$ function computations. For Leja, $\varphi_1\left(\mathcal{J}(u^n) \Delta t \right) f(u^n) \Delta t$ is computed in vertical along with the first terms of the two internal stages. The remaining terms, where the $\varphi_l(z)$ functions are applied to linear combination of nonlinear remainders, are computed individually, i.e. five $\varphi_l(z)$ function computations. 

\newcommand{\Drawx}[3]{%
  \begin{tikzpicture}[overlay,remember picture]
    \draw[#3] ([yshift = 13pt, xshift = 0pt]#1.north west) rectangle ([yshift = -2pt, xshift = 22.5pt]#2.south east);
  \end{tikzpicture}
}

\newcommand{\Drawy}[3]{%
  \begin{tikzpicture}[overlay,remember picture]
    \draw[#3] ([yshift = 13pt, xshift = 2pt]#1.north west) rectangle ([yshift = -6pt, xshift = -2pt]#2.south east);
  \end{tikzpicture}
}

\newcommand{\Drawz}[3]{%
  \begin{tikzpicture}[overlay,remember picture]
    \draw[#3] ([yshift = 13pt, xshift = 2pt]#1.north west) rectangle ([yshift = -6pt, xshift = 0pt]#2.south east);
  \end{tikzpicture}
}

\begin{align}
	a^n = u^n & + \, \tikzmark{Begin_1} \tikzmark{Begin_2} \frac{1}{2} \varphi_1\left(\frac{1}{2} \mathcal{J}(u^n) \Delta t \right) f(u^n) \Delta t \nonumber \\
	b^n = u^n & + \, \frac{9}{10} \varphi_1\left(\frac{9}{10} \mathcal{J}(u^n) \Delta t \right) f(u^n) \Delta t \tikzmark{End_1} + \tikzmark{Begin_3} \left(\frac{27}{25} \varphi_3\left(\frac{1}{2}\mathcal{J}(u^n) \Delta t\right) + \frac{729}{125} \varphi_3\left(\frac{9}{10} \mathcal{J}(u^n) \Delta t\right) \right) \mathcal{R}(a^n) \Delta t \tikzmark{End_3} \nonumber \\
	u_3^{n + 1} = u^n & + \tikzmark{Begin_4} \, \varphi_1\left(\mathcal{J}(u^n) \Delta t\right) f(u^n) \Delta t \tikzmark{End_2} \qquad \, \, + \tikzmark{Begin_6} \varphi_3(\mathcal{J}(u^n) \Delta t) \left(2\mathcal{R}(a^n) + \frac{150}{81}\mathcal{R}(b^n) \right) \Delta t  \tikzmark{End_4} \tikzmark{End_6} \\
	u_5^{n + 1} = u^n & + \tikzmark{Begin_5} \varphi_1\left(\mathcal{J}(u^n) \Delta t\right) f(u^n) \Delta t + \tikzmark{Begin_7} \varphi_3(\mathcal{J}(u^n) \Delta t)\left(18 \mathcal{R}(a^n) - \frac{250}{81} \mathcal{R}(b^n)\right) \Delta t \tikzmark{End_7} \nonumber \\
	& \qquad \qquad \qquad \qquad \qquad + \tikzmark{Begin_8} \varphi_4(\mathcal{J}(u^n) \Delta t)\left(-60 \mathcal{R}(a^n) + \frac{500}{27} \mathcal{R}(b^n)\right) \Delta t \tikzmark{End_8} \tikzmark{End_5}
	\label{eq:exprb5s3}
	\Drawa{Begin_1}{End_1}{red}
    \Drawx{Begin_2}{End_2}{blue}
    \Drawa{Begin_3}{End_3}{blue}
    \Drawy{Begin_4}{End_4}{red}
    \Drawy{Begin_5}{End_5}{red}
    \Drawz{Begin_6}{End_6}{blue}
    \Drawz{Begin_7}{End_7}{blue}
    \Drawz{Begin_8}{End_8}{blue}
\end{align}

\begin{figure}
    \centering
	\includegraphics[width = \textwidth]{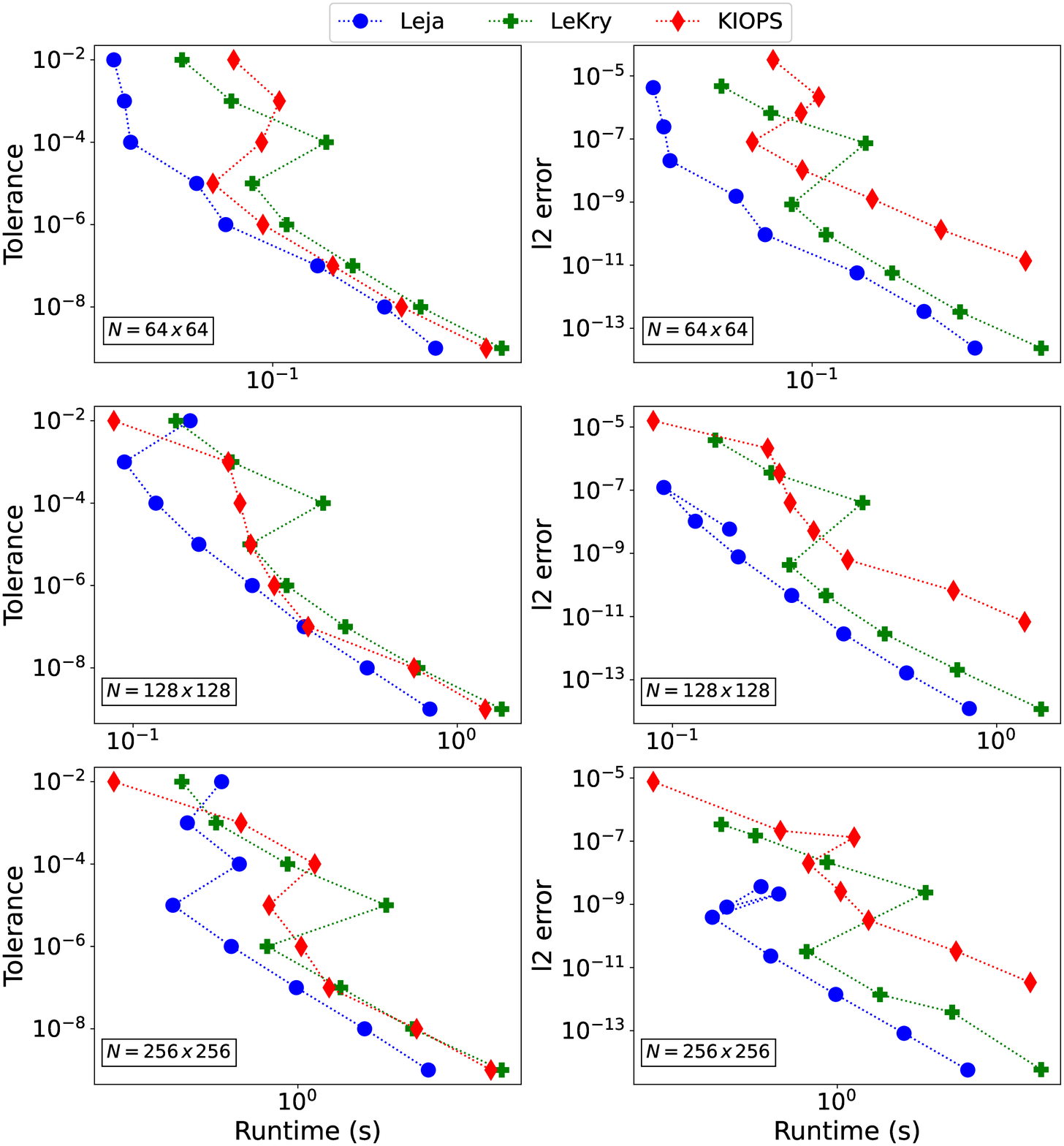}
    \caption{Comparison of the runtimes for a given value of tolerance (left panel) and the l2 norm of the global error incurred (right panel) for Leja (blue circles), LeKry (green pluses), and \kiops (red diamonds) with the EXPRB53s3 integrator for the ADR equation ($\alpha = 0.01$).}
    \label{fig:adr_53s3_0.01}
\end{figure}

The simulation results with this integrator are illustrated in Figs. \ref{fig:adr_53s3_0.01}, \ref{fig:adr_53s3_0.1}, \ref{fig:bru_53s3_0.1}, \ref{fig:sl_53s3}, and Fig. 33 - 40 in the supplementary material. For the ADR equation, Leja outperforms \kiops by a substantial margin both in terms of cost and accuracy of the solution. Here, LeKry is unable to improve upon the performance of Leja. In the case of the Allen-Cahn and Gray-Scott equations (Figs. 33 - 35 and 38 - 40), the Leja scheme has the shortest runtimes whilst ensuring that the global error incurred is less than the user-defined tolerance. Furthermore, Leja shows a relatively linear dependence of runtime as a function of the tolerance (and error). For low and intermediate stiffness (Brusselator equation, $\alpha = 0.001 - 0.01$,  Figs. 36 - 37, supplementary material), Leja and LeKry have similar runtimes for a given value of the user-defined tolerance, and both these schemes shown better accuracy than \kiops. For highly stiff cases (Fig. \ref{fig:bru_53s3_0.1}), all iterative schemes incur similar computational cost with LeKry yielding the most accurate solutions amongst the three. In all of these cases, the performance of Leja is better than or equal to than that of \kiops even though it needs an additional computation of $\varphi_l(z)$ function. Finally, LeKry is reasonably cheaper than \kiops for the semilinear equation, whilst maintaining the desired accuracy of the solution (Fig. \ref{fig:sl_53s3}). Similar to most of the previous cases, Leja is not competitive for this problem. All in all, to a reasonably good approximation, one can consider LeKry to the preferred iterative scheme for EXPRB53s3.

\begin{figure}
    \centering
	\includegraphics[width = \textwidth]{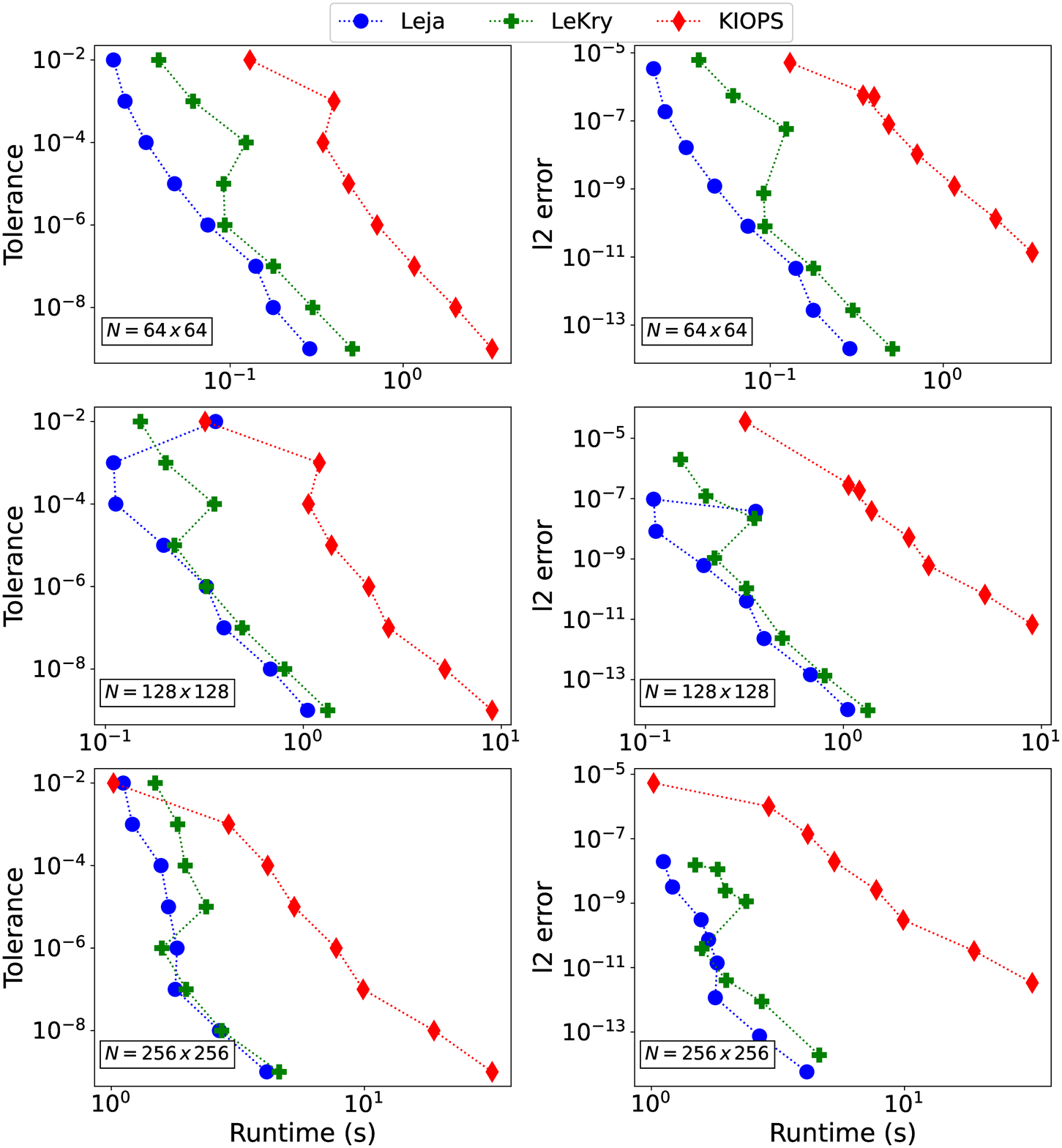} 
    \caption{Comparison of the runtimes for a given value of tolerance (left panel) and the l2 norm of the global error incurred (right panel) for Leja (blue circles), LeKry (green pluses), and \kiops (red diamonds) with the EXPRB53s3 integrator for the ADR equation ($\alpha = 0.1$).}
    \label{fig:adr_53s3_0.1}
\end{figure}

\begin{figure}
    \centering
	\includegraphics[width = \textwidth]{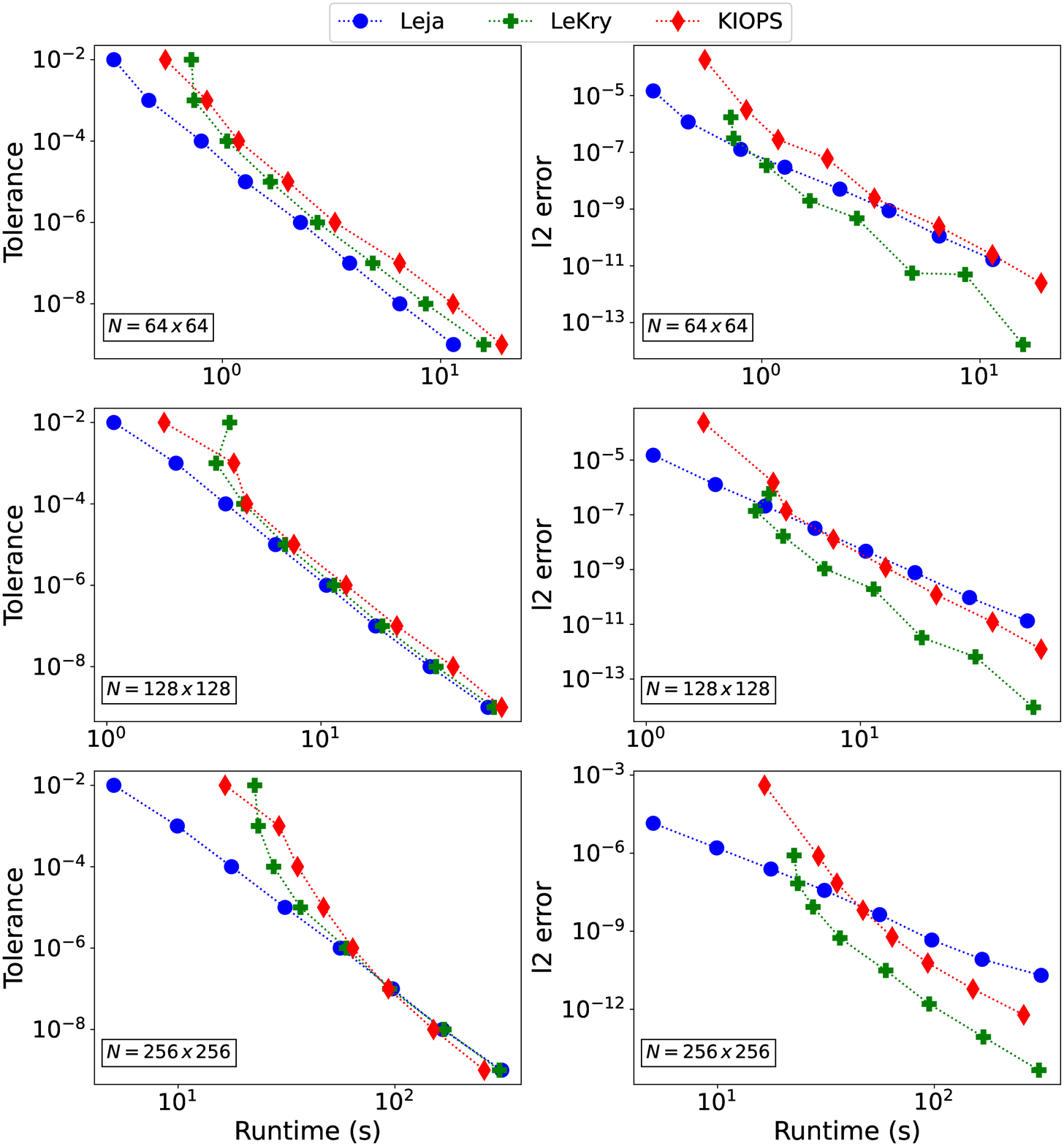} 
    \caption{Comparison of the runtimes for a given value of tolerance (left panel) and the l2 norm of the global error incurred (right panel) for Leja (blue circles), LeKry (green pluses), and \kiops (red diamonds) with the EXPRB53s3 integrator for the Brusselator equation ($\alpha = 0.1$).}
    \label{fig:bru_53s3_0.1}
\end{figure}

\begin{figure}
    \centering
	\includegraphics[width = \textwidth]{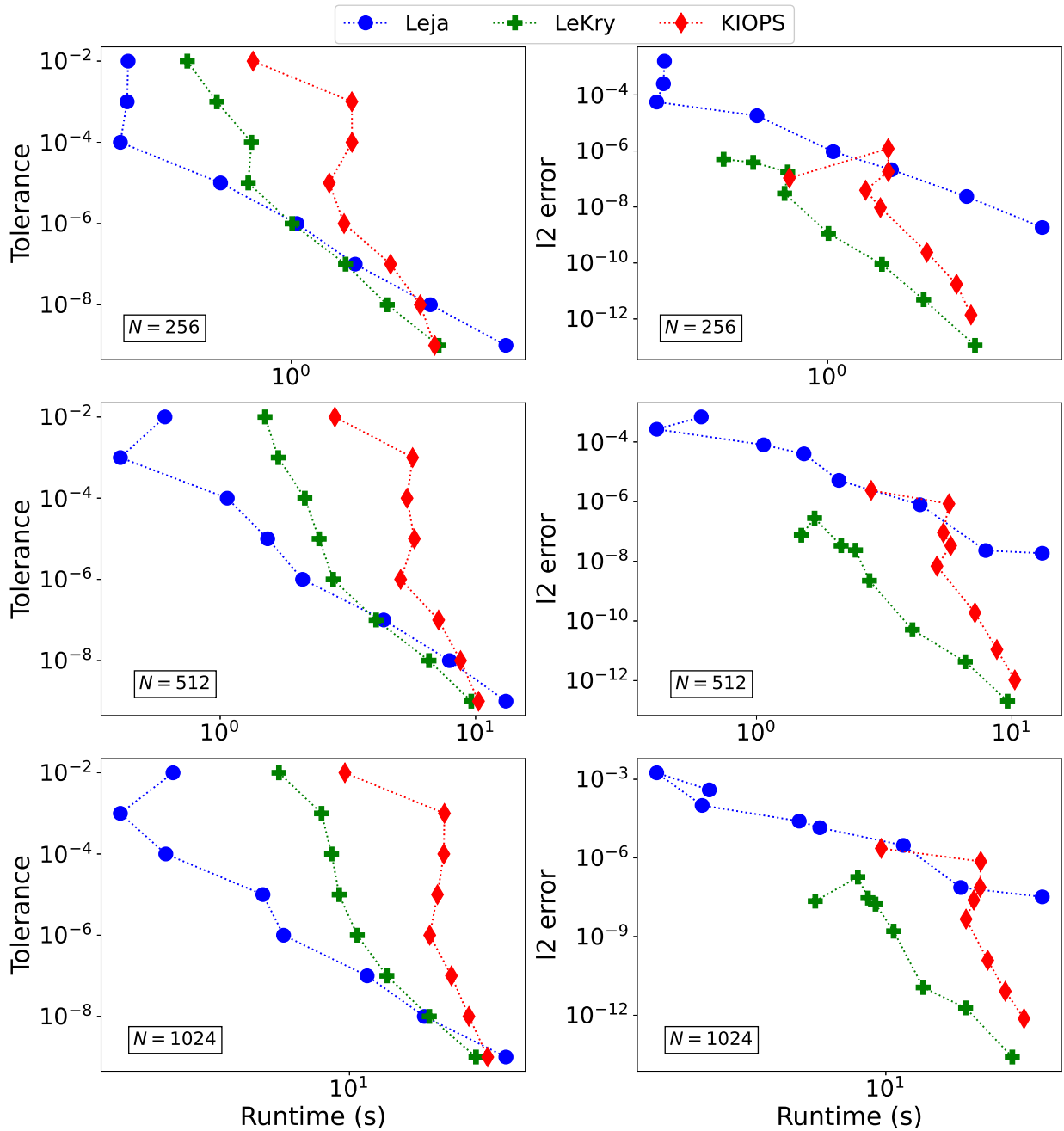} 
    \caption{Comparison of the runtimes for a given value of tolerance (left panel) and the l2 norm of the global error incurred (right panel) for Leja (blue circles), LeKry (green pluses), and \kiops (red diamonds) with the EXPRB53s3 integrator for the semilinear equation.}
    \label{fig:sl_53s3}
\end{figure}


\section{Comparison of Exponential Integrators}
\label{sec:compare_expint}

\begin{figure}[t]
    \centering
	\includegraphics[width = 0.875\textwidth]{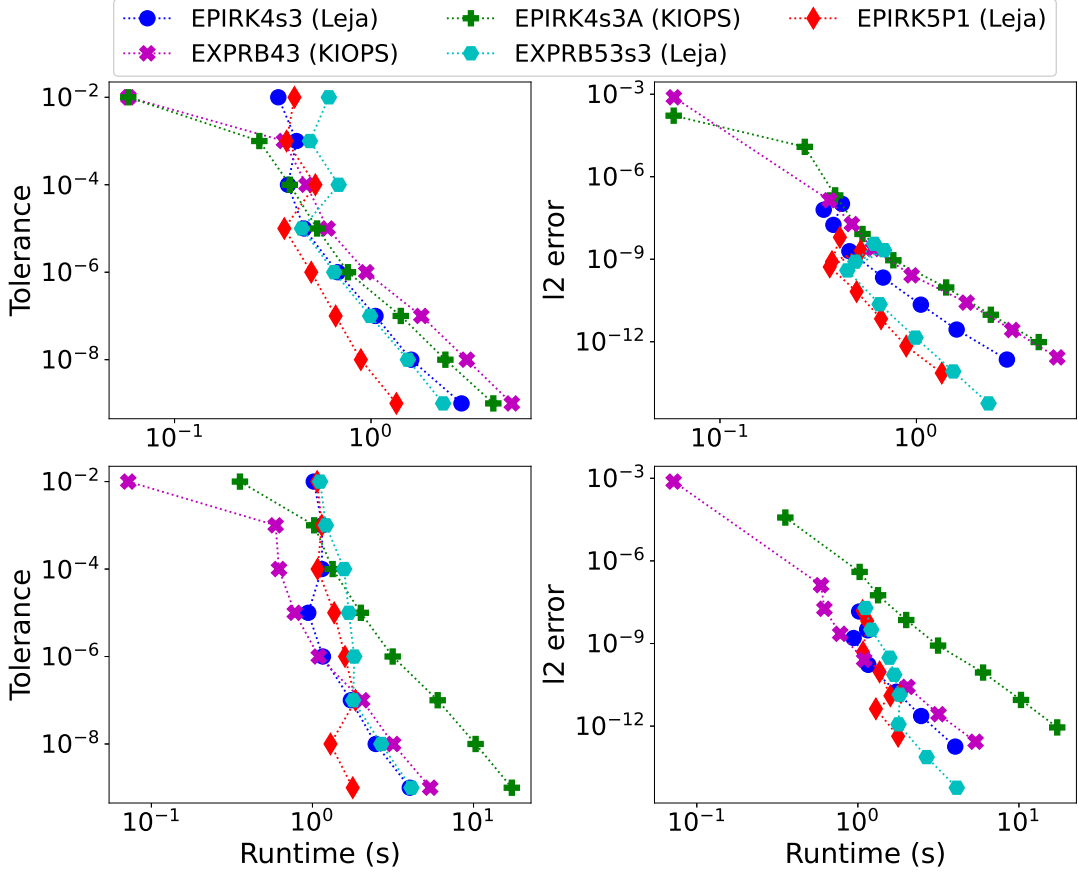} 
    \caption{Comparison of the runtimes for a given value of tolerance (left panel) and the l2 norm of the global error incurred (right panel) for the different exponential integrators for the ADR equation ($\alpha = 0.01$ (top row), $\alpha = 0.1$ (bottom row), and $N = 256 \times 256$).}
    \label{fig:adr_expint}
\end{figure}

We compare the performance of the different exponential integrators, considered in this study, for the relevant iterative scheme. We define the relevant iterative scheme as the scheme that has shown the best performance for a given integrator. We limit ourselves to Leja and \kiops for most of these comparisons as these are the two fundamental iterative schemes; LeKry is merely a combination of the two. Based on the results obtained in the previous section, we use the Leja scheme for EPIRK4s3, EPIRK5P1, and EXPRB53s3 and \kiops for EPIRK4s3A and EXPRB43. In Fig. \ref{fig:adr_expint}, we compare and contrast the runtimes of the different integrators for the ADR equation ($\alpha = 0.01$ for the top panel and $\alpha = 0.1$ for the bottom panel). It is clear that for any given value of the user-defined tolerance, the error incurred by all of the integrators is less than the tolerance. So, to compare the performance, one can simply analyse the runtimes of the integrators. In the low stiffness case, the fifth-order integrators, i.e., EPIRK5P1 and EXPRB53s3 have the shortest runtimes. EPIRK4s3 shows runtimes similar to that of EXPRB53s3 whilst incurring a slightly larger error. This is to be expected as fourth-order integrators are expected to be less accurate than fifth-order ones. The \kiops-based EPIRK4s3A and EXPRB43 are the most expensive ones. For larger values of diffusion coefficient, we see almost similar results for all integrators but EXPRB43, which has noticeably shorter runtimes than EPIRK4s3A. 

\begin{figure}[t]
    \centering
	\includegraphics[width = 0.875\textwidth]{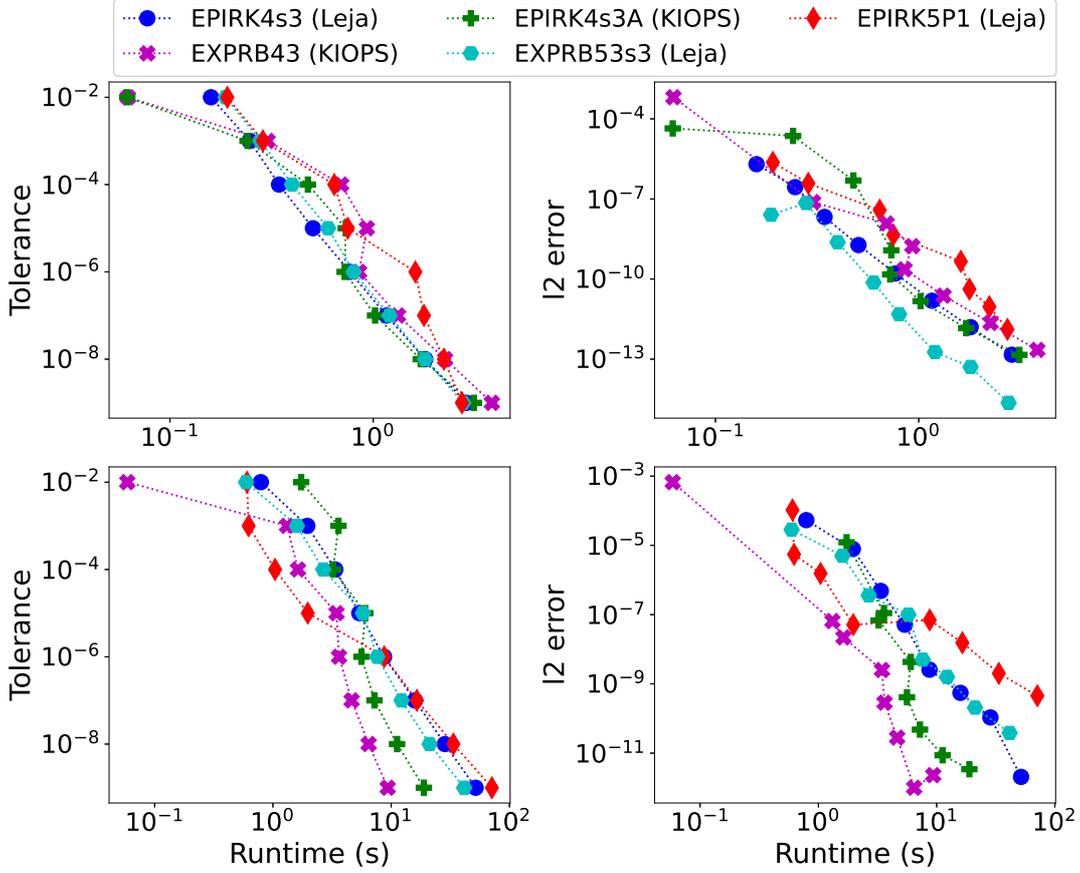} 
    \caption{Comparison of the runtimes for a given value of tolerance (left panel) and the l2 norm of the global error incurred (right panel) for the different exponential integrators for the Allen--Cahn equation ($\alpha = 0.001$ (top row), $\alpha = 0.1$ (bottom row), and $N = 256 \times 256$).}
    \label{fig:ac_expint}
\end{figure}

\begin{figure}[t]
    \centering
	\includegraphics[width = 0.875\textwidth]{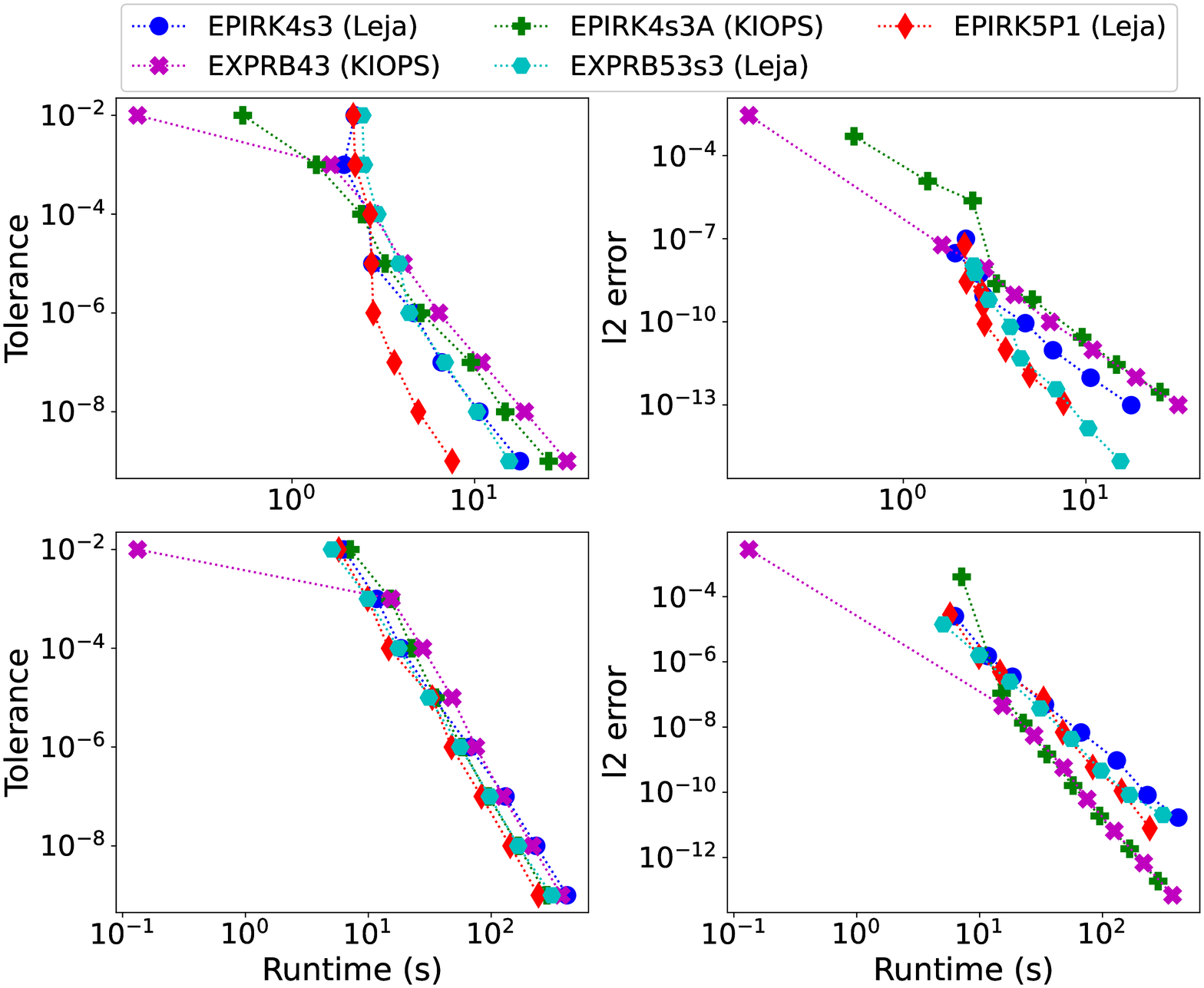} 
    \caption{Comparison of the runtimes for a given value of tolerance (left panel) and the l2 norm of the global error incurred (right panel) for the different exponential integrators for the Brusselator equation ($\alpha = 0.001$ (top row), $\alpha = 0.1$ (bottom row), and $N = 256 \times 256$).}
    \label{fig:bru_expint}
\end{figure}

\begin{figure}[t]
    \centering
	\includegraphics[width = 0.875\textwidth]{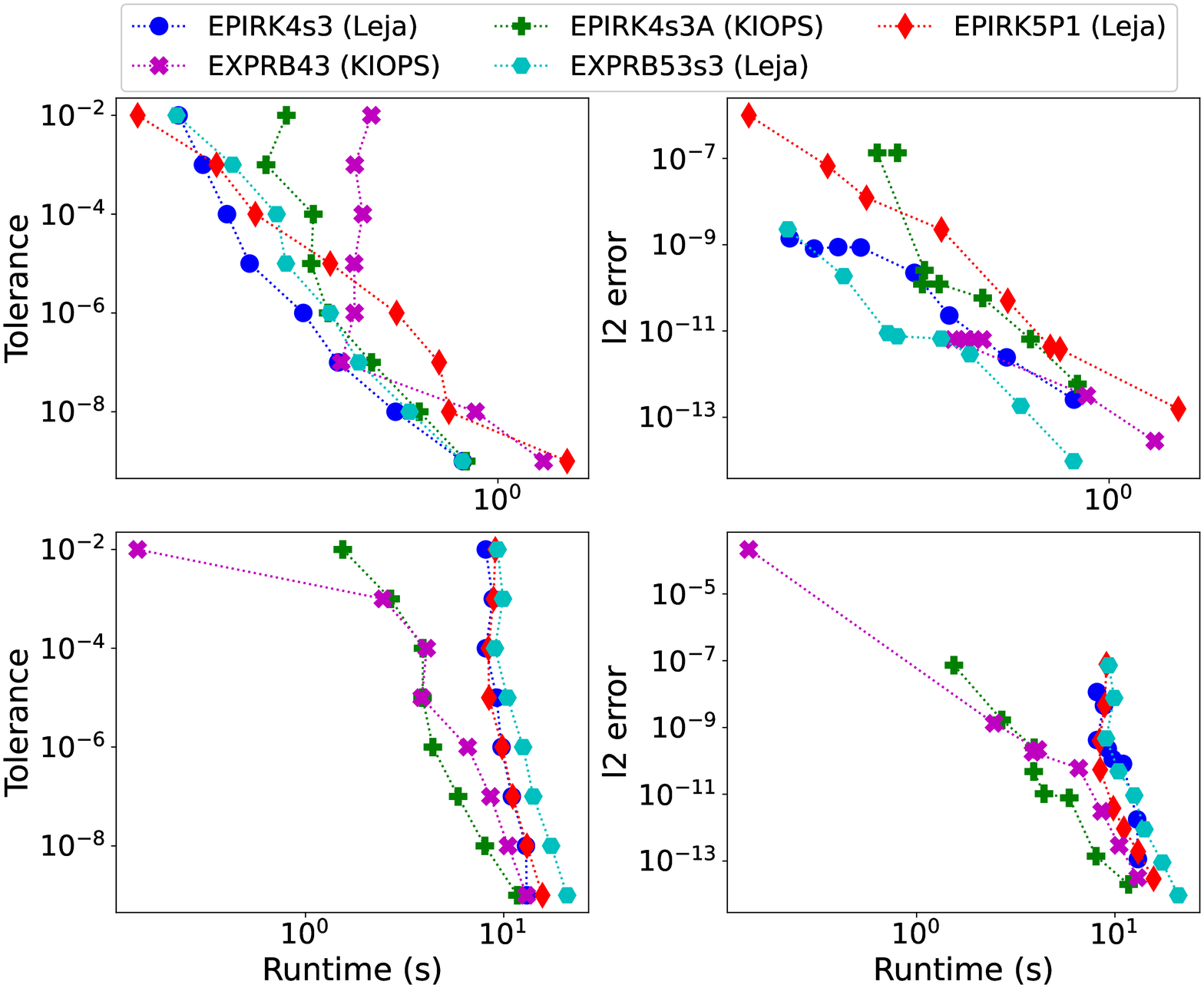}
    \caption{Comparison of the runtimes for a given value of tolerance (left panel) and the l2 norm of the global error incurred (right panel) for the different exponential integrators for the Gray--Scott equation ($\alpha = 0.001$ (top row), $\alpha = 0.1$ (bottom row), and $N = 256 \times 256$).}
    \label{fig:gs_expint}
\end{figure}

In the case of the Allen-Cahn equation (Fig. \ref{fig:ac_expint}, $\alpha = 0.001$ for the top panel and $\alpha = 0.1$ for the bottom panel), we see marginal differences in the runtimes of the different integrators, for low amounts of diffusion, with EXPRB53s3 yielding the most accurate solutions at similar expenses. With the increase in the diffusion coefficient, we see that EXPRB43 is clearly the best integrator, followed by EPIRK4s3A. They outperform the Leja based integrators. EPIRK5P1 is the most expensive amongst all of them and this is likely because of a larger number of step size rejections. Similar trends can be seen for the Brusselator and the Gray-Scott equations (Figs. \ref{fig:bru_expint} and \ref{fig:gs_expint}, respectively). Leja-based integrators are more efficient in cases with small diffusion whereas the \kiops-based ones are preferable for problems where diffusion is large. 

In the case of Gray-Scott, Bursselator, and Allen-Cahn equations, we find that the \kiops-based integrator have better performance than the Leja-based ones, whereas we find the opposite results in the case of the ADR equations. We note that the differences in the runtimes amongst the various integrators are (usually) relatively moderate. It is expected that any of these integrators are clearly expected to surpass the commonly used implicit and explicit integrators for most, if not all of the cases consider in this work (for example, see \citep{Cox02, Caliari04, Tokman06, Tokman11, Niesen09, Carr13, Luan17}). This is true especially for cases where an excellent preconditioner is either not known or extremely tedious to derive.


\section{Conclusions}
\label{sec:conclude}

Motivated by our previous work on the performance comparison of Leja and Krylov-based algorithms to compute the exponential-like functions appearing in exponential integrators, we aimed to identify the superior iterative scheme for two classes of exponential integrators: EPIRK and EXPRB. We compared the performance of the state-of-the-art \kiops algorithm, significantly different from the one used in \citet{Deka22a}, with that of the Leja interpolation scheme for a number of test problems. \kiops computes a linear combination of $\varphi_l(z)$ functions as the exponential of an augmented matrix, with effective substepping to ensure convergence and it does not require any knowledge of the spectrum of the underlying matrix prior to conducting the simulations. The Leja method, however, does need some estimate of the spectrum of the underlying matrix and we evaluate the individual $\varphi_l(z)$ functions, for enhanced computational performance.

For a performance comparison of \kiops and Leja, we considered a range of problems with varying amounts of stiffness and accuracy for several EPIRK and EXPRB methods. We did not find conclusive evidence of any one of the iterative schemes performing consistently better over the other for any specific class of exponential integrators, which renders us unable to generalise a preferred iterative scheme for a certain class of exponential integrators. However, to a reasonable degree of certainty, we have been able to identify a preferable iterative scheme for each for the individual EPIRK and EXPRB integrators that we consider in this work. We also find that the performance of a certain iterative scheme or an integrator depends on the problem and its parameters under consideration.

For the first time, we tried to combine the Leja and \kiops schemes for an exponential integrator by computing the internal stages, in vertical, using the Leja method whilst computing the linear combination of the $\varphi_l(z)$ functions using \kiops. We find remarkably excellent results for a couple of integrators, namely EPIRK4s3 and EXPRB53s3. In the future, we aim to take this a step further by implementing the Leja scheme for the small Hessenberg matrix that is obtained during the orthogonalisation process in Krylov-based methods. This is worth investigating in the massively diffusion dominated cases or scenarios with varying amounts of diffusion along different directions i.e. anisotropic diffusion.


\section*{Software}
The code used to perform the simulations with the Leja interpolation algorithm, in \texttt{Matlab}, has been added as an extension to the Krylov-based \texttt{EPIC} library. The \texttt{Python} version of the Leja method for exponential integrators is a publicly available software. Further details on the algorithm and the code can be found in Deka et al. \cite{Deka22_lexint}.

\section*{Acknowledgements}
This work has been supported by the Austrian Science Fund (FWF) project id: P32143-N32. We would like to thank Marco Caliari for providing us with the code to compute Leja points. PJD would like to thank Valentin Dallerit and Stephane Gaudreault for the useful discussions and for providing the necessary scripts of \texttt{EPIC}. This work was carried out during PJD's visit to Mayya Tokman's group at the University of California, Merced.


\appendix

\section{Computing linear combinations of \texorpdfstring{$\varphi_l$} {TEXT} %
        functions for Leja interpolation}
\label{app:linear_Leja_phi}

\begin{figure}[t]
    \centering
	\includegraphics[width = 0.875\textwidth]{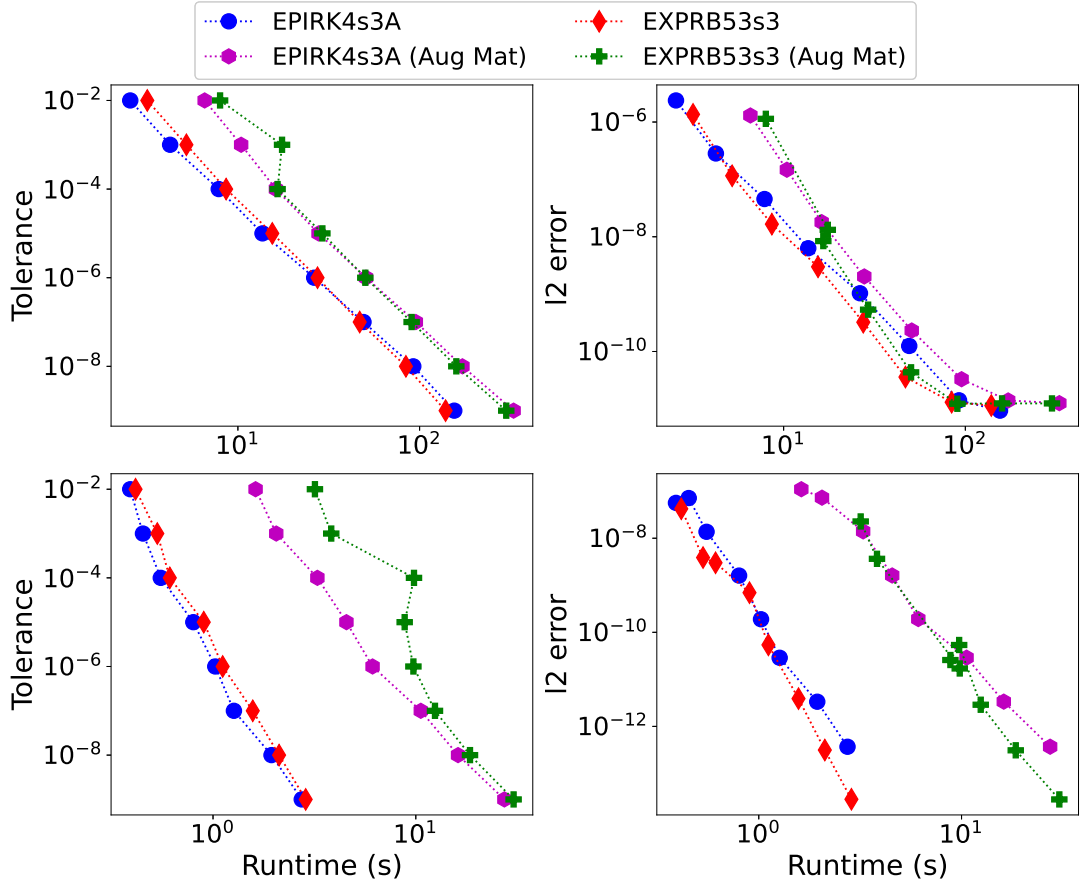}
    \caption{We compare the performance of EPIRK4s3A whilst considering an augmented matrix (magenta hexagons) with that of computing the $\varphi_l(z)$ functions individually (blue circles) for the Brusselator (top row) and Gray--Scott (bottom row) equations ($\alpha = 0.01, N = 256 \times 256$). Computing the exponential an augmented matrix clearly results in performance loss. Similar results can be seen for EXPRB53s3 where the red diamonds correspond to the case where the $\varphi_l(z)$ functions are computed separately, and the green pluses correspond to computing the augmented matrix.}
    \label{fig:aug_mat}
\end{figure}

\kiops uses a linear combination of $\varphi_l(z)$ to evaluate the various stages of an exponential integrator. However, using the same algorithm for the Leja interpolation method leads to a loss in computational performance. We show this for EPIRK4s3A and EXPRB53s3 where we have computed the third-order solution for EPIRK4s3A and the third- and fifth-order solutions for EXPRB53s3 as a linear combination of the $\varphi_l(z)$ functions, where we compute the exponential of the augmented matrix. We show, in Fig. \ref{fig:aug_mat}, that computing the exponential of the augmented matrix results in an increase in the computational runtime. Substepping of a given step size only results in an increase in the cost incurred; the accuracy requirements are already met. It is well known that the convergence of $\varphi_l(\mathcal{J} \Delta t) b_k$ depends on the magnitude of $l$, i.e., the order of the $\varphi$ function and the magnitude of $b_k$, along  with the spectrum of $\mathcal{J} \Delta t$. Higher order $\varphi$ functions have a faster convergence rate than lower-order ones. Additionally, smaller the magnitude of $b_k$, faster the convergence. The combined effect of higher order $\varphi$ functions applied to the differences of nonlinear remainders, the case of exponential integrators, result in extremely fast convergence of these functions. The most expensive part of an exponential integrator is the computation of $\varphi_1(\mathcal{J} \Delta t) f(u) \Delta t$, and as already stated $\varphi_1$ is likely to converge faster than $\varphi_0$ (the matrix exponential). This is why we obtain not-so-ideal results when we compute the exponential of the augmented matrix.


\bibliographystyle{elsarticle-harv} 
\bibliography{ref}

\begin{thebibliography}{35}
\expandafter\ifx\csname natexlab\endcsname\relax\def\natexlab#1{#1}\fi
\providecommand{\url}[1]{\texttt{#1}}
\providecommand{\href}[2]{#2}
\providecommand{\path}[1]{#1}
\providecommand{\DOIprefix}{doi:}
\providecommand{\ArXivprefix}{arXiv:}
\providecommand{\URLprefix}{URL: }
\providecommand{\Pubmedprefix}{pmid:}
\providecommand{\doi}[1]{\href{http://dx.doi.org/#1}{\path{#1}}}
\providecommand{\Pubmed}[1]{\href{pmid:#1}{\path{#1}}}
\providecommand{\bibinfo}[2]{#2}
\ifx\xfnm\relax \def\xfnm[#1]{\unskip,\space#1}\fi
\bibitem[{Al-Mohy and Higham(2011)}]{Higham10}
\bibinfo{author}{Al-Mohy, A.H.}, \bibinfo{author}{Higham, N.J.},
  \bibinfo{year}{2011}.
\newblock \bibinfo{title}{Computing the action of the matrix exponential, with
  an application to exponential integrators}.
\newblock \bibinfo{journal}{SIAM J. Sci. Comput.} \bibinfo{volume}{33},
  \bibinfo{pages}{488--511}.
\newblock \DOIprefix\doi{10.1137/100788860}.
\bibitem[{Arnoldi(1951)}]{Arnoldi1951}
\bibinfo{author}{Arnoldi, W.E.}, \bibinfo{year}{1951}.
\newblock \bibinfo{title}{The principle of minimized iterations in the solution
  of the matrix eigenvalue problem}.
\newblock \bibinfo{journal}{Q. Appl. Math.} \bibinfo{volume}{9},
  \bibinfo{pages}{17--29}.
\newblock \DOIprefix\doi{10.1090/qam/42792}.
\bibitem[{Baglama et~al.(1998)Baglama, Calvetti and Reichel}]{Baglama98}
\bibinfo{author}{Baglama, J.}, \bibinfo{author}{Calvetti, D.},
  \bibinfo{author}{Reichel, L.}, \bibinfo{year}{1998}.
\newblock \bibinfo{title}{Fast leja points}.
\newblock \bibinfo{journal}{Electron. Trans. Numer. Anal.} \bibinfo{volume}{7},
  \bibinfo{pages}{124 -- 140}.
\newblock \URLprefix \url{http://eudml.org/doc/119747}.
\bibitem[{Bates et~al.(2009)Bates, Brown and Han}]{Bates09}
\bibinfo{author}{Bates, P.W.}, \bibinfo{author}{Brown, S.},
  \bibinfo{author}{Han, J.}, \bibinfo{year}{2009}.
\newblock \bibinfo{title}{{Numerical Analysis for a Nonlocal Allen-Cahn
  Equation}}.
\newblock \bibinfo{journal}{Int. J. Numer. Anal. Model.} \bibinfo{volume}{6},
  \bibinfo{pages}{33--49}.
\newblock \URLprefix
  \url{http://global-sci.org/intro/article_detail/ijnam/755.html}.
\bibitem[{Bergamaschi et~al.(2006)Bergamaschi, Caliari, Martinez and
  Vianello}]{Bergamaschi06}
\bibinfo{author}{Bergamaschi, L.}, \bibinfo{author}{Caliari, M.},
  \bibinfo{author}{Martinez, A.}, \bibinfo{author}{Vianello, M.},
  \bibinfo{year}{2006}.
\newblock \bibinfo{title}{Comparing leja and krylov approximations of large
  scale matrix exponentials}.
\newblock \bibinfo{journal}{Proc. ICCS} ,
  \bibinfo{pages}{685--692}\DOIprefix\doi{10.1007/11758549_93}.
\bibitem[{Caliari and Ostermann(2009)}]{Caliari09}
\bibinfo{author}{Caliari, M.}, \bibinfo{author}{Ostermann, A.},
  \bibinfo{year}{2009}.
\newblock \bibinfo{title}{{Implementation of exponential Rosenbrock-type
  integrators}}.
\newblock \bibinfo{journal}{Appl. Numer. Math.} \bibinfo{volume}{59},
  \bibinfo{pages}{568 -- 581}.
\newblock \DOIprefix\doi{10.1016/j.apnum.2008.03.021}.
\bibitem[{Caliari et~al.(2004)Caliari, Vianello and Bergamaschi}]{Caliari04}
\bibinfo{author}{Caliari, M.}, \bibinfo{author}{Vianello, M.},
  \bibinfo{author}{Bergamaschi, L.}, \bibinfo{year}{2004}.
\newblock \bibinfo{title}{Interpolating discrete advection–diffusion
  propagators at leja sequences}.
\newblock \bibinfo{journal}{J. Comput. Appl. Math.} \bibinfo{volume}{172},
  \bibinfo{pages}{79 -- 99}.
\newblock \DOIprefix\doi{10.1016/j.cam.2003.11.015}.
\bibitem[{{Caliari} et~al.(2007){Caliari}, {Vianello} and
  {Bergamaschi}}]{Caliari07b}
\bibinfo{author}{{Caliari}, M.}, \bibinfo{author}{{Vianello}, M.},
  \bibinfo{author}{{Bergamaschi}, L.}, \bibinfo{year}{2007}.
\newblock \bibinfo{title}{{The LEM exponential integrator for
  advection-diffusion-reaction equations}}.
\newblock \bibinfo{journal}{J. Comput. Appl. Math.} \bibinfo{volume}{210},
  \bibinfo{pages}{56--63}.
\newblock \DOIprefix\doi{10.1016/j.cam.2006.10.055}.
\bibitem[{Carr et~al.(2013)Carr, Moroney and Turner}]{Carr13}
\bibinfo{author}{Carr, E.J.}, \bibinfo{author}{Moroney, T.J.},
  \bibinfo{author}{Turner, I.W.}, \bibinfo{year}{2013}.
\newblock \bibinfo{title}{Performance assessment of exponential rosenbrock
  methods for large systems of ode}, in: \bibinfo{booktitle}{Proceedings of the
  16th Biennial Computational Techniques and Applications Conference,
  CTAC-2012}, pp. \bibinfo{pages}{C102--C118}.
\bibitem[{Cox and Matthews(2002)}]{Cox02}
\bibinfo{author}{Cox, S.}, \bibinfo{author}{Matthews, P.},
  \bibinfo{year}{2002}.
\newblock \bibinfo{title}{Exponential time differencing for stiff systems}.
\newblock \bibinfo{journal}{J. Sci. Comput.} \bibinfo{volume}{176},
  \bibinfo{pages}{430 -- 455}.
\newblock \DOIprefix\doi{10.1006/jcph.2002.6995}.
\bibitem[{Deka and Einkemmer(2022a)}]{Deka22a}
\bibinfo{author}{Deka, P.J.}, \bibinfo{author}{Einkemmer, L.},
  \bibinfo{year}{2022}a.
\newblock \bibinfo{title}{{Efficient adaptive step size control for exponential
  integrators}}.
\newblock \bibinfo{journal}{Comput. Math. Appl.} \bibinfo{volume}{123},
  \bibinfo{pages}{59--74}.
\newblock \DOIprefix\doi{10.1016/j.camwa.2022.07.011}.
\bibitem[{Deka and Einkemmer(2022b)}]{Deka22b}
\bibinfo{author}{Deka, P.J.}, \bibinfo{author}{Einkemmer, L.},
  \bibinfo{year}{2022}b.
\newblock \bibinfo{title}{{Exponential} {Integrators} for {Resistive}
  {Magnetohydrodynamics}: {Matrix-free} {Leja} {Interpolation} and {Efficient}
  {Adaptive} {Time} {Stepping}}.
\newblock \bibinfo{journal}{ApJS} \bibinfo{volume}{259}, \bibinfo{pages}{57}.
\newblock \DOIprefix\doi{10.3847/1538-4365/ac5177}.
\bibitem[{{Deka} et~al.(2022){Deka}, {Einkemmer} and {Tokman}}]{Deka22_lexint}
\bibinfo{author}{{Deka}, P.J.}, \bibinfo{author}{{Einkemmer}, L.},
  \bibinfo{author}{{Tokman}, M.}, \bibinfo{year}{2022}.
\newblock \bibinfo{title}{{LeXInt: Package for Exponential Integrators
  employing Leja interpolation}}.
\newblock \bibinfo{journal}{arXiv e-prints} ,
  \bibinfo{pages}{arXiv:2208.08269}\href{http://arxiv.org/abs/2208.08269}{{\tt
  arXiv:2208.08269}}.
\bibitem[{Einkemmer(2018)}]{Einkemmer18}
\bibinfo{author}{Einkemmer, L.}, \bibinfo{year}{2018}.
\newblock \bibinfo{title}{An adaptive step size controller for iterative
  implicit methods}.
\newblock \bibinfo{journal}{Appl. Numer. Math.} \bibinfo{volume}{132},
  \bibinfo{pages}{182 -- 204}.
\newblock \DOIprefix\doi{10.1016/j.apnum.2018.06.002}.
\bibitem[{{Gaudreault} et~al.(2018){Gaudreault}, {Rainwater} and
  {Tokman}}]{Gaudreault18}
\bibinfo{author}{{Gaudreault}, S.}, \bibinfo{author}{{Rainwater}, G.},
  \bibinfo{author}{{Tokman}, M.}, \bibinfo{year}{2018}.
\newblock \bibinfo{title}{{KIOPS: A fast adaptive Krylov subspace solver for
  exponential integrators}}.
\newblock \bibinfo{journal}{JCP} \bibinfo{volume}{372},
  \bibinfo{pages}{236--255}.
\newblock \DOIprefix\doi{10.1016/j.jcp.2018.06.026},
  \href{http://arxiv.org/abs/1804.05126}{{\tt arXiv:1804.05126}}.
\bibitem[{Gray and Scott(1984)}]{Gray84}
\bibinfo{author}{Gray, P.}, \bibinfo{author}{Scott, S.}, \bibinfo{year}{1984}.
\newblock \bibinfo{title}{{Autocatalytic reactions in the isothermal,
  continuous stirred tank reactor: Oscillations and instabilities in the system
  A + 2B → 3B; B → C}}.
\newblock \bibinfo{journal}{Chem. Eng. Sci.} \bibinfo{volume}{39},
  \bibinfo{pages}{1087--1097}.
\newblock \DOIprefix\doi{10.1016/0009-2509(84)87017-7}.
\bibitem[{Hairer and Wanner(1996)}]{HairerII}
\bibinfo{author}{Hairer, E.}, \bibinfo{author}{Wanner, G.},
  \bibinfo{year}{1996}.
\newblock \bibinfo{title}{Solving Ordinary Differential Equations II, Stiff and
  Differential-Algebraic Problems}.
\newblock \bibinfo{publisher}{Springer Berlin Heidelberg}.
\newblock \URLprefix
  \url{https://link.springer.com/book/10.1007/978-3-642-05221-7}.
\bibitem[{Hochbruck and Ostermann(2005)}]{Hochbruck05}
\bibinfo{author}{Hochbruck, M.}, \bibinfo{author}{Ostermann, A.},
  \bibinfo{year}{2005}.
\newblock \bibinfo{title}{Explicit exponential runge--kutta methods for
  semilinear parabolic problems}.
\newblock \bibinfo{journal}{SIAM J. Numer. Anal.} \bibinfo{volume}{43},
  \bibinfo{pages}{1069--1090}.
\newblock \DOIprefix\doi{10.1137/040611434}.
\bibitem[{Hochbruck and Ostermann(2010)}]{Ostermann10}
\bibinfo{author}{Hochbruck, M.}, \bibinfo{author}{Ostermann, A.},
  \bibinfo{year}{2010}.
\newblock \bibinfo{title}{Exponential integrators}.
\newblock \bibinfo{journal}{Acta Numer.} \bibinfo{volume}{19},
  \bibinfo{pages}{209 -- 286}.
\newblock \DOIprefix\doi{10.1017/S0962492910000048}.
\bibitem[{Hochbruck et~al.(2009)Hochbruck, Ostermann and
  Schweitzer}]{Hochbruck09}
\bibinfo{author}{Hochbruck, M.}, \bibinfo{author}{Ostermann, A.},
  \bibinfo{author}{Schweitzer, J.}, \bibinfo{year}{2009}.
\newblock \bibinfo{title}{Exponential rosenbrock-type methods}.
\newblock \bibinfo{journal}{SIAM J. Numer. Anal.} \bibinfo{volume}{47},
  \bibinfo{pages}{786--803}.
\newblock \DOIprefix\doi{10.1137/080717717}.
\bibitem[{{Lefever} and {Nicolis}(1971)}]{Lefever71}
\bibinfo{author}{{Lefever}, R.}, \bibinfo{author}{{Nicolis}, G.},
  \bibinfo{year}{1971}.
\newblock \bibinfo{title}{{Chemical instabilities and sustained oscillations}}.
\newblock \bibinfo{journal}{J. Theor. Biol.} \bibinfo{volume}{30},
  \bibinfo{pages}{267--284}.
\newblock \DOIprefix\doi{10.1016/0022-5193(71)90054-3}.
\bibitem[{Leja(1957)}]{Leja1957}
\bibinfo{author}{Leja, F.}, \bibinfo{year}{1957}.
\newblock \bibinfo{title}{Sur certaines suites liées aux ensembles plans et
  leur application à la représentation conforme}.
\newblock \bibinfo{journal}{Ann. Polon. Math.} \bibinfo{volume}{4},
  \bibinfo{pages}{8--13}.
\newblock \URLprefix \url{http://eudml.org/doc/208291}.
\bibitem[{Loffeld and Tokman(2013)}]{Tokman13}
\bibinfo{author}{Loffeld, J.}, \bibinfo{author}{Tokman, M.},
  \bibinfo{year}{2013}.
\newblock \bibinfo{title}{{Comparative performance of exponential, implicit,
  and explicit integrators for stiff systems of ODEs}}.
\newblock \bibinfo{journal}{J. Comput. Appl. Math.} \bibinfo{volume}{241},
  \bibinfo{pages}{45--67}.
\newblock \DOIprefix\doi{10.1016/j.cam.2012.09.038}.
\bibitem[{Luan(2017)}]{Luan17}
\bibinfo{author}{Luan, V.T.}, \bibinfo{year}{2017}.
\newblock \bibinfo{title}{{Fourth-order two-stage explicit exponential
  integrators for time-dependent PDEs}}.
\newblock \bibinfo{journal}{Appl. Numer. Math.} \bibinfo{volume}{112},
  \bibinfo{pages}{91}.
\newblock \DOIprefix\doi{10.1016/j.apnum.2016.10.008}.
\bibitem[{Luan and Ostermann(2014)}]{Luan14}
\bibinfo{author}{Luan, V.T.}, \bibinfo{author}{Ostermann, A.},
  \bibinfo{year}{2014}.
\newblock \bibinfo{title}{{Exponential Rosenbrock methods of order five —
  construction, analysis and numerical comparisons}}.
\newblock \bibinfo{journal}{J. Comput. Appl. Math.} \bibinfo{volume}{255},
  \bibinfo{pages}{417--431}.
\newblock \DOIprefix\doi{10.1016/j.cam.2013.04.041}.
\bibitem[{Michels et~al.(2017)Michels, Luan and Tokman}]{Tokman17a}
\bibinfo{author}{Michels, D.L.}, \bibinfo{author}{Luan, V.T.},
  \bibinfo{author}{Tokman, M.}, \bibinfo{year}{2017}.
\newblock \bibinfo{title}{A stiffly accurate integrator for elastodynamic
  problems}.
\newblock \bibinfo{journal}{ACM Trans. Graph.} \bibinfo{volume}{36}.
\newblock \DOIprefix\doi{10.1145/3072959.3073706}.
\bibitem[{Moler and Van~Loan(2003)}]{Moler03}
\bibinfo{author}{Moler, C.}, \bibinfo{author}{Van~Loan, C.},
  \bibinfo{year}{2003}.
\newblock \bibinfo{title}{Nineteen dubious ways to compute the exponential of a
  matrix, twenty-five years later}.
\newblock \bibinfo{journal}{SIAM Rev.} \bibinfo{volume}{45},
  \bibinfo{pages}{3--49}.
\newblock \DOIprefix\doi{10.1137/S00361445024180}.
\bibitem[{Moler and Loan(1978)}]{Moler78}
\bibinfo{author}{Moler, C.B.}, \bibinfo{author}{Loan, C.V.},
  \bibinfo{year}{1978}.
\newblock \bibinfo{title}{Nineteen dubious ways to compute the exponential of a
  matrix}.
\newblock \bibinfo{journal}{Siam Rev.} \bibinfo{volume}{20},
  \bibinfo{pages}{801--836}.
\newblock \DOIprefix\doi{10.1137/1020098}.
\bibitem[{Niesen and Wright(2012)}]{Niesen09}
\bibinfo{author}{Niesen, J.}, \bibinfo{author}{Wright, W.M.},
  \bibinfo{year}{2012}.
\newblock \bibinfo{title}{Algorithm 919: A krylov subspace algorithm for
  evaluating the $\varphi$-functions appearing in exponential integrators}.
\newblock \bibinfo{journal}{ACM Trans. Math. Softw.} \bibinfo{volume}{38}.
\newblock \DOIprefix\doi{10.1145/2168773.2168781}.
\bibitem[{{Rainwater} and {Tokman}(2016)}]{Tokman16}
\bibinfo{author}{{Rainwater}, G.}, \bibinfo{author}{{Tokman}, M.},
  \bibinfo{year}{2016}.
\newblock \bibinfo{title}{{A new approach to constructing efficient stiffly
  accurate EPIRK methods}}.
\newblock \bibinfo{journal}{JCP} \bibinfo{volume}{323},
  \bibinfo{pages}{283--309}.
\newblock \DOIprefix\doi{10.1016/j.jcp.2016.07.026}.
\bibitem[{{Rainwater} and {Tokman}(2017)}]{Tokman17b}
\bibinfo{author}{{Rainwater}, G.}, \bibinfo{author}{{Tokman}, M.},
  \bibinfo{year}{2017}.
\newblock \bibinfo{title}{{Designing efficient exponential integrators with
  EPIRK framework}}, in: \bibinfo{booktitle}{International Conference of
  Numerical Analysis and Applied Mathematics (ICNAAM 2016)}, p.
  \bibinfo{pages}{020007}.
\newblock \DOIprefix\doi{10.1063/1.4992153}.
\bibitem[{Tokman(2006)}]{Tokman06}
\bibinfo{author}{Tokman, M.}, \bibinfo{year}{2006}.
\newblock \bibinfo{title}{{Efficient integration of large stiff systems of ODEs
  with exponential propagation iterative (EPI) methods}}.
\newblock \bibinfo{journal}{JCP} \bibinfo{volume}{213},
  \bibinfo{pages}{748--776}.
\newblock \DOIprefix\doi{10.1016/j.jcp.2005.08.032}.
\bibitem[{Tokman(2011)}]{Tokman11}
\bibinfo{author}{Tokman, M.}, \bibinfo{year}{2011}.
\newblock \bibinfo{title}{{A new class of exponential propagation iterative
  methods of Runge–Kutta type (EPIRK)}}.
\newblock \bibinfo{journal}{JCP} \bibinfo{volume}{230},
  \bibinfo{pages}{8762--8778}.
\newblock \DOIprefix\doi{10.1016/j.jcp.2011.08.023}.
\bibitem[{Tokman et~al.(2012)Tokman, Loffeld and Tranquilli}]{Tokman12}
\bibinfo{author}{Tokman, M.}, \bibinfo{author}{Loffeld, J.},
  \bibinfo{author}{Tranquilli, P.}, \bibinfo{year}{2012}.
\newblock \bibinfo{title}{{New Adaptive Exponential Propagation Iterative
  Methods of Runge--Kutta Type}}.
\newblock \bibinfo{journal}{SIAM J. Sci. Comput.} \bibinfo{volume}{34},
  \bibinfo{pages}{A2650--A2669}.
\newblock \DOIprefix\doi{10.1137/110849961}.
\bibitem[{{Van Der Vorst}(1987)}]{Vorst87}
\bibinfo{author}{{Van Der Vorst}, H.}, \bibinfo{year}{1987}.
\newblock \bibinfo{title}{An iterative solution method for solving f(a)x = b,
  using krylov subspace information obtained for the symmetric positive
  definite matrix a}.
\newblock \bibinfo{journal}{J. Comput. Appl. Math.} \bibinfo{volume}{18},
  \bibinfo{pages}{249--263}.
\newblock \DOIprefix\doi{10.1016/0377-0427(87)90020-3}.

\end{thebibliography}

\newpage





\end{document}